\documentclass[12pt]{amsart}
\usepackage[margin=1.2in]{geometry}
\usepackage{color}

\title[Affine Dynamical Sampling on Graphs]
{Randomized Space-Time Sampling for Affine Graph Dynamical Systems}

\author[L. Gong, L. Huang]{Le Gong, Longxiu Huang}

\usepackage{hyperref}
\usepackage{cleveref}

\date{\vspace{-5ex}}
\usepackage{amsthm, amsfonts, mathrsfs}
\usepackage{ dsfont }
\usepackage{amssymb}
\usepackage{enumerate}
\usepackage{graphicx}
\usepackage{subfig}

\usepackage{float}
\usepackage{bbm}
\usepackage{amsmath}
\usepackage{comment}
\usepackage{hyperref}
\usepackage{listings}
\usepackage{color}
\usepackage{todonotes}

\usepackage{diagbox}
\usepackage[dvipsnames]{xcolor}
\usepackage{algorithm}
\usepackage{booktabs}
\usepackage{algorithmicx}
\usepackage[noend]{algpseudocode}
\usepackage{tikz}
\allowdisplaybreaks

\hypersetup{
    colorlinks,
    citecolor=blue,
    filecolor=blue,
    linkcolor=blue,
    urlcolor=blue
}
\numberwithin{equation}{section}
\newtheorem{theorem}{Theorem}
\newtheorem{proposition}[theorem]{Proposition}
\newtheorem{lemma}[theorem]{Lemma}

\newtheorem{definition}{Definition}
\newtheorem{remark}{Remark}

\newcommand{\R}{\mathbb{R}}

\newcommand{\N}{\mathbb{N}}

\newcommand{\bx}{\textbf{x}}
\newcommand{\by}{\textbf{y}}
\newcommand{\bz}{\textbf{z}}
\newcommand{\bw}{\textbf{w}}
\newcommand{\be}{\textbf{e}}
\newcommand{\bv}{\textbf{v}}

\usepackage[foot]{amsaddr}

\definecolor{aacolor}{rgb}{0.05, 0.75, 1}

\begin{document}
	\date{}

\address{\textrm{(Le Gong)}
	Department of Mathematics,
	Vanderbilt University,
	Nashville, TN 37240 USA}
\email{le.gong@vanderbilt.edu}

\address{\textrm{(Longxiu Huang)}
        Department of Computational Mathematics, Science and Engineering \& Department of
        Mathematics, Michigan State University,
        East Lansing, MI 48824, USA}
\email{huangl3@msu.edu}

\keywords{Bandlimited graph signals; Affine dynamical system; Source recovery; Random   Space-time sampling;  Compressed sensing. }
\subjclass [2010] {05C50, 94A12, 94A20}

\maketitle
\begin{abstract}
 This paper investigates the problem of dynamical sampling for graph signals influenced by a constant source term. We consider signals evolving over time according to a linear dynamical system on a graph, where both the initial state and the source term are bandlimited. We introduce two random space-time sampling regimes and analyze the conditions under which stable recovery is achievable. While our framework extends recent work on homogeneous dynamics, it addresses a fundamentally different setting where the evolution includes a constant source term. This results in a non-orthogonal-diagonalizable  system matrix, rendering classical spectral techniques inapplicable and introducing new challenges in sampling design, stability analysis, and joint recovery of both the initial state and the forcing term. A key component of our analysis is the spectral graph weighted coherence, which characterizes the interplay between the sampling distribution and the graph structure. We establish sampling complexity bounds ensuring stable recovery via the Restricted Isometry Property (RIP), and develop a robust recovery algorithm with provable error guarantees. The effectiveness of our method is validated through extensive experiments on both synthetic and real-world datasets.
\end{abstract}

\section{Introduction}
Dynamical sampling is a fundamental problem in signal processing that seeks to recover an evolving signal from a limited number of space-time samples. Unlike traditional sampling theory, which focuses on static signals, dynamical sampling leverages the temporal evolution of signals to enable efficient reconstruction. Given a signal \( x_t \) residing in a  Hilbert space and evolving under a  linear dynamical system

\begin{equation}
    \bx_t = A \bx_{t-1},
\end{equation}  
where \( A \) is a known evolution operator, the goal is to determine  space-time sampling strategies that allow for stable and accurate reconstruction of the initial signal.  

The concept of  dynamical sampling  was introduced by  Aldroubi et al.  \cite{aldroubi2017dynamical,aldroubi2013dynamical}, motivated by the pioneering work of  Lu et al. \cite{lu2009spatial}, which studied space-time sampling problems in an evolutionary system. Aldroubi et al. established a mathematical framework to characterize feasible sampling sets using tools from frame theory and operator theory in harmonic analysis. Further developments and extensions to different settings have been explored in \cite{ dokmanic2011sensor,murray2015estimating, murray2016physics, ranieri2011sampling}.  

Most existing research \cite{aceska2013dynamical, aldroubi2021sampling,aldroubi2019frames,aldroubi2020phaseless, cabrelli2017dynamical, christensen2019frame,  lai2017undersampled,  tang2017universal, ulanovskii2021reconstruction} has focused on deterministic sampling frameworks, where the spatial sampling locations remain fixed over time. These studies typically consider signals evolving over well-structured domains, such as Euclidean spaces or separable Hilbert spaces. However, in many modern applications—such as social networks, sensor grids, or biological systems—signals are naturally supported on graphs, where the domain is determined by an underlying network topology.

To address this, recent work \cite{huang2025random,huang2024robust,yao2023space} investigates the problem of sampling and reconstructing dynamical graph signals that are sparse and/or bandlimited. These works propose some random sampling regimes and introduce the notion of spectral graph weighted coherence to quantify the sufficient sampling size required for stable reconstruction from partial observations. While these models effectively capture the structural properties of graph signals, they assume homogeneous evolution, i.e., signal propagation driven solely by the system dynamics.

However, in many real-world scenarios, the signal evolution is also influenced by external source terms, which may represent environmental inputs, control actions, or spontaneous signal generation. To capture such effects, it is essential to incorporate a source term into the dynamical model, leading to the extended formulation \[\bx_{t+1}=A\bx_{t}+\bw_t\] where 
$\bw_t$  denotes the source input at time $t$. Modeling the system in this way enables a more realistic representation of evolving signals in practical applications \cite{aastrom2021feedback,lian2021robustness}. Thus, extending dynamical sampling theory to incorporate source terms is critical for broader applicability, and it introduces new theoretical and algorithmic challenges related to identifiability and reconstruction stability.
In this work, we focus on a specific setting where the source term is constant over time, i.e.,  
\begin{equation}\label{model:discrete_source}
    \mathbf{x}_{t+1} = A \mathbf{x}_t + \mathbf{w},
\end{equation}
which has been recently studied in \cite{aldroubi2024dynamical, aldroubi2024characterization}. We assume that both the initial signal \( \mathbf{x}_0 \) and the source term \( \mathbf{w} \) are bandlimited signals over a graph, and that the driving operator \( A \) is known in advance. The signal evolution is observed up to \( s \) time instances. We define the space-time sampling set as  
\[
\Omega = \bigcup_{t=0}^{s-1} \left( \Omega_t \times \{t\} \right) \subseteq [n] \times \{0, 1, \dots, s-1\},
\]  
where \( \Omega_t = \{\omega_{1,t}, \omega_{2,t}, \dots, \omega_{m_t,t}\} \subseteq [n] \) denotes the set of observed spatial locations at time \( t \).  We propose and study two random space-time sampling regimes:  
\begin{enumerate}
    \item[(i)] \textbf{Fixed spatial sampling across time:} the spatial sampling set remains constant, i.e., \( \Omega_t = \Omega_0 \) for all \( t \in \{0, \dots, s-1\} \);
    \item[(ii)] \textbf{Independent random sampling at each time:} the spatial sampling set \( \Omega_t \) is drawn independently at each time step.
\end{enumerate}

Our objective is to analyze the conditions under which stable recovery of both the bandlimited initial state and the constant source term is achievable under randomized space-time sampling regimes. We aim to establish rigorous theoretical guarantees that characterize the interplay between the signal dynamics, the forcing term, and the sampling strategy. Although our formulation for the affine system with a source term bears superficial resemblance to the homogeneous setting considered in \cite{huang2024robust}, the extension is far from trivial. In the homogeneous case, the system admits a clean spectral analysis through repeated application of a diagonalizable operator \( A \). In contrast, the affine system can be rewritten as a block evolution:
\[
\begin{bmatrix}
x_t \\
w
\end{bmatrix}
= 
\begin{bmatrix}
A & I \\
0 & I
\end{bmatrix}^t
\begin{bmatrix}
x_0 \\
w
\end{bmatrix},
\]
where the block matrix is, in general, \emph{not orthogonally diagonalizable}. This structural change invalidates the spectral decomposition techniques and RIP analysis tools employed in \cite{huang2024robust}. As a result, new mathematical techniques are required to analyze stability and sampling complexity. Furthermore, we tackle the more challenging task of jointly recovering both \( x_0 \) and \( w \), which introduces a bilinear inverse problem structure with increased ill-posedness and sensitivity to noise, particularly under random measurements.
\subsection{Related work}
\subsubsection{Sampling for bandlimited graph signals}

Sampling and reconstruction methods for bandlimited graph signals have significantly progressed, typically categorized into two main groups: selection sampling and aggregation sampling. 

 Selection sampling includes deterministic methods (e.g.,~\cite{chen2015discrete,chamon2017greedy,pesenson2008sampling,anis2016efficient,huang2020reconstruction}) and random methods inspired by compressed sensing (e.g.,~\cite{chen2015signal,hashemi2018accelerated,puy2018random}). However, these approaches primarily focus on static graph signals. Recently, Huang et al.~\cite{huang2024robust} generalized Vandergheynst et al.'s random sampling methods~\cite{puy2018random} to bandlimited diffusion fields through three novel random space-time sampling schemes. Aggregation sampling reconstructs signals using linear combinations from neighboring nodes based on a graph shift operator that mirrors the adjacency structure. The deterministic scheme proposed in~\cite{marques2015sampling}, along with its extensions in~\cite{wang2016local, yang2021orthogonal}, can be viewed as  our sampling regime 1, where the sampling locations remain fixed across all time instances. In contrast, the randomized approach in~\cite{valsesia2018sampling} employs a graph shift operator with i.i.d.~Gaussian entries.

A detailed survey can be found in Tanaka~\cite{tanaka2020sampling}. In this work, we  extend the first two random sampling strategies from~\cite{huang2024robust} to scenarios involving bandlimited diffusion fields with an added fixed source term.

\subsubsection{Dynamical sampling for  structured signals}


Dynamical sampling addresses the problem of recovering an evolving signal from spatiotemporal samples. When the underlying signal exhibits structure---such as bandlimitedness, shift-invariance, or a graph-based domain---this framework enables efficient sampling strategies that leverage both temporal dynamics and structural priors.

For bandlimited signals, Aldroubi et al.~\cite{aldroubi2021sampling} and Al-Hammali and Faridani~\cite{al2023uniform} investigate recovery from spatiotemporal measurements, with an emphasis on stability and convergence. Their results show that nonuniform or augmented sampling schemes—such as periodic perturbations or additional spatial samples—can enable stable reconstruction and yield rapidly converging series, even under sub-Nyquist sampling or in the presence of blind spots.
 In the case of shift-invariant spaces, which arise naturally in time-frequency analysis and image processing, the works \cite{aceska2014dynamical,huo2022multivariate,zhang2017dynamical} extend classical sampling theory by incorporating dynamical evolution. They provide recovery guarantees and constructive algorithms for both univariate and multivariate settings. This framework has also been extended to high-dimensional signals. For instance, Wang et.al.~\cite{wang2025three} develops a dynamical sampling theory for three-dimensional signal recovery, where the evolution is modeled by tensor t-product operators. The work provides necessary conditions on the sampling set and formulates the reconstruction problem as an efficiently solvable optimization task, validated by numerical experiments.

Dynamical sampling has also been generalized to signals defined on graphs, where the spatial structure is discrete and irregular. The works \cite{aldroubi2024reconstructing,huang2025random,huang2024robust} propose robust recovery strategies for graph signals using spatiotemporal samples, drawing on spectral graph theory and compressed sensing techniques. These methods are well-suited for applications such as sensor networks and transportation systems, where data acquisition is limited and the underlying geometry is non-Euclidean. Building on these developments, this work extends randomized sampling regimes from the recovery of bandlimited graph signals—serving as initial conditions of linear dynamical systems—to the more general task of reconstructing both the initial state and a constant source term in affine linear systems.
\subsubsection{Dynamical sampling with source term}
Several recent works have extended the dynamical sampling framework to systems driven by external source terms. Aldroubi et al.~\cite{aldroubi2023predictive} address the recovery of burst-like forcing terms in an initial value problem (IVP). They introduce two specialized classes of samplers that enable accurate prediction of the system's evolution during intervals without bursts and develop robust algorithms for recovering the burst-like input, even in the presence of noise and a large background signal. Related works~\cite{aldroubi2023recovery,aldroubi2024reconstruction} further examine the recovery of rapidly decaying source terms in IVP settings.

Recently, Cheng et al.~\cite{cheng2021estimate} consider  Model~\ref{model:discrete_source}, focusing on the case where the evolution operator $A$ is unknown. Their work aims to determine sampling conditions under which the spectrum of $A$ can be recovered from the observed measurements.

In a related direction, Aldroubi et al.~\cite{aldroubi2024dynamical, aldroubi2024characterization} present a comprehensive analysis of recovering the constant source term in Model~\ref{model:discrete_source}, considering both finite- and infinite-dimensional separable Hilbert spaces. Using spatiotemporal measurements, they derive necessary and sufficient conditions on the observations to ensure exact recovery. Building on this line of work, Vashisht et al.~\cite{vashisht2025frames} recently examine the stability of recovering source terms in discrete dynamical systems indexed over non-uniform discrete sets, which arise from spectral pairs in infinite-dimensional separable Hilbert spaces.

In this work, we further consider the problem of simultaneously recovering both the initial state and a constant source term from randomized spatiotemporal samples, under the assumption that both the initial state and the source term are bandlimited. We also investigate the stability of the proposed reconstruction method.

\subsection{Summary of contributions} 
This work introduces a new framework for dynamical sampling of graph signals influenced by a constant source term. In contrast to prior work that considers purely homogeneous evolution, we study the recovery of both the initial state and a constant source term from random space-time samples. Our contributions can be summarized as follows.
\begin{itemize}
    \item Incorporating source dynamics into random dynamical sampling: 
    We extend randomized dynamical sampling to graph signals governed by affine dynamics with a constant source term. Unlike the homogeneous model in \cite{huang2024robust}, our framework handles a non-orthogonal-diagonalizable  system and requires a fundamentally new stability analysis.  Note that in the case where $\|A\|<1$,  the homogeneous dynamics studied in~\cite{huang2024robust} exhibit a decay in information from the initial state over time. In contrast, our model reveals a trade-off between this decay and the accumulation of information from the source term. 
    
    \item Simultaneous recovery of initial state and source term: 
    In contrast to \cite{aldroubi2024dynamical, aldroubi2024characterization}, which provide deterministic recovery guarantees under known sampling patterns, our approach supports randomized sampling and simultaneously recovers both the initial state \( \bx_0 \) and the source term \( \bw \). Moreover, while their methods do not account for measurement noise, our proposed algorithm explicitly incorporates noise and provides corresponding error analysis.
    
    \item Two random space-time sampling regimes and theoretical guarantees: 
    We propose and analyze two random sampling strategies: (i) fixed spatial sampling across time, and (ii) independent random sampling at each time instance. For both regimes, we introduce the concept of \emph{spectral graph weighted coherence} to characterize the interaction between the signal structure and sampling distribution. Using this quantity, we establish sampling complexity bounds under which the RIP holds with high probability, ensuring stable recovery. The details are presented in \Cref{secrip}. 
    
    \item Robust reconstruction algorithm and error analysis: 
    We develop a practical reconstruction algorithm based on \( \ell_2 \)-regularization that leverages both the graph Laplacian and the bandlimited structure of the signals. Our method provides provable guarantees under noisy measurements, and we analyze the reconstruction error in terms of the noise level and spectral properties of the graph.
 \item Extensive experiments on synthetic and real datasets: 
    We validate our method through numerical experiments on both synthetic and real-world datasets. The results demonstrate exact recovery in the noiseless case and  robustness under noise. In real data settings, our approach outperforms existing methods designed for bandlimited signals, supporting the practical impact of our theoretical contributions. More details are illustrated in~\Cref{experiments}.
\end{itemize}


\section{Preliminaries, problem setting and notation}
\subsection{Graph operations} We consider a simple, undirected, weighted graph $\mathcal{G} = (\mathcal{V}, \mathcal{E})$, where the node set is $\mathcal{V} = \{v_1, v_2, \dots, v_n\}$ and the edge set is $\mathcal{E}\subset\mathcal{V}\times \mathcal{V}$.  Let \(\{v_i, v_j\}\) denote an edge connecting vertices \(v_i\) and \(v_j\). Since the graph is undirected, we do not distinguish between \((v_i, v_j)\) and \((v_j, v_i)\), and use \(\{v_i, v_j\}\) to represent the edge. The weighted adjacency matrix \(W \in \mathbb{R}^{n \times n}\) is defined as:
\[
W(i, j) \triangleq 
\begin{cases}
\alpha_{ij}, & \text{if } \{v_i, v_j\} \in \mathcal{E}, \\
0, & \text{otherwise},
\end{cases}
\quad \text{where } \alpha_{ij} \in \mathbb{R}_{+}, \text{ for all } v_i, v_j \in \mathcal{V}.
\]

The weight \(\alpha_{ij}\) associated with each edge typically reflects the similarity or dependency between the corresponding vertices. These weights and the overall graph connectivity are often determined by the physics or structure of the application—for example, weights may be inversely proportional to the physical distance between nodes in a network. By construction, the adjacency matrix \(W\) is symmetric.

The \textit{degree} of a vertex \(v_i \in \mathcal{V}\), denoted \(\deg(v_i)\), is defined as \(\deg(v_i) = \sum_{j=1}^n W(i,j)\). The diagonal degree matrix is given by \(D := \operatorname{diag}(\deg(v_i))_{i=1}^n\).

We now define two fundamental operators associated with the graph \(\mathcal{G}\).

\begin{definition}
The \textit{normalized diffusion operator} of a graph \(\mathcal{G}\) with weighted adjacency matrix \(W \in \mathbb{R}^{n \times n}\) is defined as
\[
N = D^{-\frac{1}{2}} W D^{-\frac{1}{2}}.
\]
We refer to \(N\) as the \textit{diffusion matrix} of \(\mathcal{G}\). The associated \textit{normalized graph Laplacian} is defined as
\[
L = I_n - N,
\]
where \(I_n\) is the \(n \times n\) identity matrix.
\end{definition}

\subsection{Bandlimited graph signals} Given a graph \(\mathcal{G} = (\mathcal{V}, \mathcal{E})\), suppose it is associated with a normalized graph Laplacian matrix \(L \in \mathbb{R}^{n \times n}\). The matrix \(L\) is typically real, symmetric, and positive semi-definite, enabling frequency analysis of signals defined on the graph. Under these conditions, \(L\) admits the eigendecomposition:
\[
L = U \Theta U^\top,
\]
where \(U \in \mathbb{R}^{n \times n}\) is an orthonormal matrix whose columns are the eigenvectors of \(L\), and \(\Theta = \text{diag}(\theta_1,\theta_2, \dots, \theta_n)\) is a diagonal matrix of nonnegative eigenvalues. Without loss of generality, we assume that the eigenvalues are ordered as \(0=\theta_1 \leq \theta_2 \leq \dots \leq \theta_n\).

This eigendecomposition provides a notion of frequency in the graph domain. The Graph Fourier Transform (GFT) of a signal \(\bx \in \mathbb{R}^n\) and its inverse are defined as:
\[
\hat{\bx} = U^\top \bx, \quad \text{and} \quad \bx = U \hat{\bx}.
\]
The diagonal entries of \(\Theta\) represent the graph frequencies, and the spectral content of \(\bx\) is captured by the vector \(\hat{\bx}\). A signal is said to be \(k\)-bandlimited if its spectral energy is concentrated entirely within the first \(k\) graph frequencies; that is, \(\hat{\bx}\) has nonzero entries only in its first \(k\) components. This leads to the following formal definition:

\begin{definition}
    A signal $\bx\in\R^n$ defined on the nodes of the graph $\mathcal{G}$ is k-bandlimited with some nonzero $k\in[n]$ if $\bx\in\text{span}(U_k),$ where $U_k:=U(:,1:k)\in\R^{n\times k}$ is the submatrix consisting of the first $k$ columns of $U.$ 
\end{definition}
For a \(k\)-bandlimited signal \(\bx \in \mathbb{R}^n\), we redefine its GFT and its inverse as
\[
\hat{\bx} := U_k^\top \bx, \text{ and }
\bx = U_k \hat{\bx},
\]
reflecting that \(\bx\) lies entirely within the subspace spanned by the first \(k\) graph Fourier modes.

\subsection{Problem setting}
We focus on a discrete, time-dependent graph signal $\bx_t \in \mathbb{R}^n$ that evolves according to the dynamics governed by an operator $A$ and a source term $\bw$. The signal evolution is described by the following equation:  
\begin{equation}
    \label{model}
    \bx_{t+1} = A \bx_t + \bw,    
\end{equation}
where $A$ determines how the signal propagates across the graph, and $\bw$ represents an external input or source term, assumed to be independent of time. Importantly, both the initial signal state $\bx_0$ and the source term $\bw$ are unknown. 

Our objective is to reconstruct $\bx_0$ and $\bw$ simultaneously in a stable way, leveraging random sampling. To facilitate the analysis, we define 
$\bw_{\text{aug}}=\begin{bmatrix}
    \bx_0\\
\bw
\end{bmatrix}$, reducing the problem to recovering 
$\bw_{\text{aug}}$. 

We observe the system dynamics over \(s\) time steps. Let \(\Omega^s = \bigcup_{t=0}^{s-1} \left( \Omega_t \times \{t\} \right) \)  denote the set of space–time sampling locations, where \(\Omega_t = \{\omega_{1,t}, \dots, \omega_{m_t,t}\} \subset [n]\) is the set of spatial observation indices at time \(t\). The sampling operator at time \(t\) is represented by the matrix \(S_t \in \mathbb{R}^{m_t \times n}\), defined as:
\[
S_t = \begin{bmatrix}
    \delta_{\omega_{1,t}},&\cdots,&
    \delta_{\omega_{m_t,t}}
\end{bmatrix}^{\top},
\]
where \(\delta_{\omega_{i,t}} \in \mathbb{R}^n\) is the Dirac delta vector (i.e., the standard basis vector) that is nonzero only at index \(\omega_{i,t}\):
\[
\delta_{\omega_{i,t}}(j) = 
\begin{cases}
1, & \text{if } j = \omega_{i,t}, \\
0, & \text{otherwise}.
\end{cases}
\]
The sampling locations are selected randomly, with replacement, according to a predefined probability distribution $\textbf{p}_t$ on $[n]$.

We observe the signal evolution process and obtain the data $\bz=\{S_t\bx_t+\be_t: t=0,1,\cdots,s-1\}$ with $\be_t$ denoting the noise at time $t$. The goal is to analyze the sampling complexities required for stable recovery of the initial signal $\bx_0$ and the source term $\bw$ from $\bz$, and to propose robust reconstruction algorithms.

In this paper, the operator \( A \in \mathbb{R}^{n \times n} \) in \eqref{model} is referred to as the \emph{graph shift operator}. We assume that \( A \) is a diagonalizable matrix sharing the same eigenspace as the graph Laplacian \( L \); typical examples include polynomial functions of \( L \), the matrix exponential \( e^{-\alpha L} \text{ with }\alpha\in\mathbb{R}\), and related matrix functions. Specifically, 
\[
A := U \Lambda U^\top,
\]
where \( \Lambda \in \mathbb{R}^{n \times n} \) is a diagonal matrix of eigenvalues with entries \( \lambda_1, \lambda_2, \ldots, \lambda_n \).

In \eqref{model}, if we assume that both the initial state $\bx_0$ and the source term $\bw$ are k-bandlimited signals, i.e., $\bx_0, \bw\in\text{span}(U_k),$ we can now claim the following proposition. The proof follows directly from the relation given in \eqref{eqn:rel} below.

\begin{proposition}
The graph signal $\bx_t$ defined by the model in \eqref{model} belongs to $\text{span}(U_k)$ for each $t\in\N$ if both the initial state $\bx_0$ and the source term $\bw$ are k-bandlimited signals.
\end{proposition}

\subsection{Notation}
In this work, we let \( s \) denote the maximum number of sampling time steps. Throughout, we use capital letters (e.g., \( A \)) to denote matrices, bold lowercase letters (e.g., \( \bx \)) for vectors, and regular lowercase letters (e.g., \( k \)) for scalars. The graph Fourier transform of a signal \( \bx \) is denoted by \( \hat{\bx} \), and \( \|\bx\| \) represents the Euclidean norm of \( \bx \). The set of the first \( d \) natural numbers is denoted by \( [d] \). 

We use \( \mathbf{0} \) to denote a matrix with all entries equal to zero, and \( I \) to denote the identity matrix, with their dimensions depending on the context. For a matrix \( U \), \( \|U\|_2 \) denotes its spectral norm, and \( \text{span}(U) \) represents the vector space spanned by its columns. The notation \( \text{diag}(\bx) \) refers to a diagonal matrix whose diagonal entries are given by the vector \( \bx \). Given a diagonal matrix \( \Lambda = \text{diag}(\lambda_{1}, \lambda_{2},\dots, \lambda_{n}) \), we define \( \Lambda_k \) for \( k \leq n \) as its top-left \( k \times k \) principal submatrix, i.e., \( \Lambda_k = \Lambda(1{:}k, 1{:}k) \). 



\section{Dynamic random sampling of $k$-bandlimited signals}\label{secrip}
\subsection{Stable embedding of $k$-bandlimited signals via space–time samples} 
Recall that $s-1\ (s\ge 1)$ is the maximum sampling time point. Assume that the eigenvalue decomposition of $A$ is of the form $A=U\Lambda U^\top$. Using the relation $\bx_{t+1} = A \bx_t + \bw$, we derive the following expression: 
\begin{align*}
    \bx_{t}& 
   =\begin{bmatrix}
        I&\bf{0}
    \end{bmatrix}\begin{bmatrix}
        A&I\\
        \bf{0}&I
    \end{bmatrix}^t\begin{bmatrix}
        \bx_0\\
        \bw
    \end{bmatrix}=\ \begin{bmatrix}
        U&\bf{0}
    \end{bmatrix}
    \begin{bmatrix}
        \Lambda&I\\
        \bf{0}&I
    \end{bmatrix}^{t}\begin{bmatrix}
        U^\top&\bf{0}\\
        \bf{0}&U^\top
    \end{bmatrix}\begin{bmatrix}
        \bx_0\\
        \bw
    \end{bmatrix}.
\end{align*}
While the above form may appear structurally similar to expressions derived in prior work such as \cite{huang2024robust}, a crucial distinction lies in the matrix $
\begin{bmatrix}
A & I \\
\mathbf{0} & I
\end{bmatrix}$, 
which is, in general, \emph{not orthogonally diagonalizable}. This lack of  orthogonally diagonalizability presents a significant analytical challenge, as standard spectral decomposition techniques that are effective in the homogeneous case cannot be directly applied. Consequently, our analysis requires a fundamentally different approach and more delicate treatment to address the intricacies introduced by the affine term.

To proceed, we define the matrix sequences
 $\Lambda^0:=I\in \mathbb{R}^{n \times n}$,  $\bar{\Lambda}^0:=\mathbf{0}\in \mathbb{R}^{n \times n}$ and $\bar{\Lambda}^t:=\sum_{j=0}^{t-1}\Lambda^j\ (t\ge 1)$. Then we have
\begin{align*}
    \bx_{t}=&
    \begin{bmatrix}
        U&\bf{0}
    \end{bmatrix}\begin{bmatrix}
        \Lambda^{t}&\bar{\Lambda}^{t}\\
        \bf{0}& I
    \end{bmatrix}\begin{bmatrix}
        U^\top&\bf{0}\\
        \bf{0}&U^\top
    \end{bmatrix}\begin{bmatrix}
        \bx_0\\
        \bw
    \end{bmatrix}
    =U\begin{bmatrix}
       \Lambda^{t}&\bar{\Lambda}^{t}
   \end{bmatrix}\begin{bmatrix}
        U^\top&\bf{0}\\
        \bf{0}&U^\top
    \end{bmatrix}\begin{bmatrix}
        \bx_0\\
        \bw
    \end{bmatrix}.
\end{align*}
 Since $x_0$ and $w$ are $k$-bandlimited, we  have 
\begin{align}\label{eqn:rel}
    \bx_{t}
    =&U_k\begin{bmatrix}
       \Lambda_k^{t}&\bar{\Lambda}_k^{t}
   \end{bmatrix}\begin{bmatrix}
        U_k^\top&\bf{0}\\
        \bf{0}&U_k^\top
    \end{bmatrix}\begin{bmatrix}
        \bx_0\\
        \bw
    \end{bmatrix}.
\end{align}
The above observations inspire us to introduce a map $\pi_{A,s}:\mathbb{R}^{2n}\rightarrow\mathbb{R}^{sn}$ defined by 
\begin{equation*}
\pi_{A,s}(\bv)=\begin{bmatrix}
U_k\begin{bmatrix}
    I &\bf{0}
\end{bmatrix}\\
U_k
\begin{bmatrix}
    \Lambda_k & I
\end{bmatrix}\\
\vdots\\
    U_k 
\begin{bmatrix}
    \Lambda_k^{s-1} & \bar{\Lambda}_k^{s-1}
\end{bmatrix}
\end{bmatrix}\begin{bmatrix}
        U_k^\top&\bf{0}\\
        \bf{0}&U_k^\top
    \end{bmatrix}
\bv
\end{equation*}
and the following theorem shows that $\pi_{A,s}$ is a stable embedding.

\begin{theorem}\label{embedding}
    For any $\bx_0$, $\bw$ being $k$-bandlimited and  $\bv=\begin{bmatrix}
        \bx_0\\
        \bw
    \end{bmatrix}\in\mathbb{R}^{2n}$,  we have that
    \begin{equation}\label{boundpi}
        c_{A,k,s}\|\bv\|^2 \le \|\pi_{A,s}(\bv)\|^2\le C_{A,k,s}\|\bv\|^2
    \end{equation}
where $0<c_{A,k,s}\le C_{A,k,s}.$
\end{theorem}

\subsection{Sampling procedure}
To select a subset of nodes for sampling, we define $\mathbf{P} = [\mathbf{p}_0,\mathbf{p}_1, \dots, \mathbf{p}_{s-1}] \in \mathbb{R}^{n\times s}$ with $\mathbf{p}_t\in\mathbb{R}^n$ being a probability distribution over the set $\{1,2, \dots, n\}$, which serves as the sampling distribution. We assume $\mathbf{p}_t(i) > 0$ for all $i\in[n]$ and $\sum_{i=1}^n \mathbf{p}_t(i) = 1$.

Let $\Omega_t=\{\omega_{1,t},\omega_{2,t},\cdots,\omega_{m_t,t}\}$ be the sampling node set at time $t$. It is constructed by drawing $m_t$ indices independently with replacement from   $[n]$ according to the probability distribution $\mathbf{p}_t,$ which means that 
\[
\mathbb{P}(\omega_{j,t}=i)=\mathbf{p}_t(i),\quad\forall j\in[m_t]\ \text{and}\   i\in[n].
\]

Next let us define the sampling matrix $S_t\in\R^{m_t\times n}.$ This matrix satisfies
\[
   S_t(i,j) = \begin{cases} 
   1 & \text{if } j = \omega_{i,t} \\
   0 & \text{otherwise}.
   \end{cases}
\]
for all $i\in[m_t]$ and $j\in[n]$. In other words,
$S_t=\begin{bmatrix}
    \delta_{\omega_{1,t}}&\cdots&
    \delta_{\omega_{m_t,t}}
\end{bmatrix}^{\top}\in\R^{m_t\times n}.$ The sample $\by_t$  at time $t\in \{0,1,\cdots,s-1\}$ is collected through $\by_t=S_t\bx_t.$

For a signal $\bx_t\in\R^n$ defined on the nodes of the graph, its sampled version $\by_t\in\R^{m_t}$ is given by
\[
\by_t(j):=\bx_t({\omega_{j,t}}),\quad\forall j\in[m_t].
\]
Note that since sampling is done with replacement, the values of $\by_j$ are not necessarily distinct for different $j$. Although this sampling procedure permits each node to be selected multiple times, in practice, each chosen node can be sampled just once, with any duplicates added artificially afterward.

\subsection{Sampling regimes}

In our analysis, we consider two distinct regimes for the sampling matrices $S_t$, which determine the structure of the observations over time.

\begin{enumerate}
    \item \textbf{Fixed Sampling Matrix}:  In this regime, a time-invariant sampling matrix is used, i.e., \( S_t = S \) for all \( t \). Consequently, the set of sampled spatial locations \( \Omega_t \) and the number of observed entries \( m_t = m \) remain fixed across time steps. We also assume that the probability distributions \( \mathbf{p}_t \) are identical for all \( t \), and the sampling locations are determined solely by \( \mathbf{p}_0 \). For clarity, we denote the full distribution matrix as \( \mathbf{P}^{(1)} \), with only the first column, \( \mathbf{P}^{(1)}(:,1) \), used to generate the sampling set. 
    The observed data are then expressed as:

    \begin{equation}
        \bz = \{ S \bx_t + \be_t : t=0,1,\cdots,s-1 \},
    \end{equation}
    where each observation is taken at the same locations in every time step. This setting has constraints due to the lack of variability in the observation pattern over time.

    \item \textbf{Time-Varying Sampling Matrix}: In contrast, this regime allows the sampling matrices $S_t$ to vary with time, meaning that the observed spatial locations $\Omega_t$ are different at each time step. As a result, the number of observed entries $m_t$ may also vary over time. The observed data in this setting is given by:
    \begin{equation}
        \bz = \{ S_t \bx_t + \be_t : t=0,1,\dots,s-1 \}.
    \end{equation}
    This regime introduces additional randomness in the sampling process, which can potentially improve reconstruction performance by covering a more diverse set of locations across time.
\end{enumerate}

To pave the way for introducing the main theorem, we first discuss an important quantity called \textbf{graph spectral weighted coherence}, a related concept previously proposed in \cite{huang2024robust}.

In this work, the graph spectral weighted coherence, denoted by $ \nu^{(j)}\ (j=1,2)$, is defined across two distinct regimes. In regime 1, $\nu^{(1)} $ is determined by a summation over different time steps, while in regime 2, $ \nu^{(2)} $ captures the coherence at each time step, both of which depend on the probability distributions $\mathbf{P}^{(1)}$ and $\mathbf{P}^{(2)}$, respectively. These coherence values guide us in determining the required number of samples at each time step to ensure accurate signal reconstruction.


\begin{definition}[Graph spectral weighted coherence]\label{coherence} Let
    $\mathbf{P}^{(j)}$ be the sampling distribution for $j=1,2$. The graph spectral weighted coherence $\nu^{(j)}$ is defined as follows:
    \begin{enumerate}
        \item Regime 1:
        $\nu^{(1)}=
        \max\limits_{1\le \ell \le n}
    \left\{\frac{1}{\mathbf{p}_0^{(1)}(\ell)}\left\|\sum\limits_{t=0}^{s-1}\begin{bmatrix}
    \Lambda_k^{t} \\
    \bar{\Lambda}_k^{t}
\end{bmatrix}
U_k^{\top}\delta_{\ell}\delta_{\ell}^{\top}U_k
\begin{bmatrix}
    \Lambda_k^{t} & \bar{\Lambda}_k^{t}\\
\end{bmatrix}\right\|_2\right\}.$
    \item 
    Regime 2: $\nu^{(2)}=\begin{bmatrix}
        \nu^{(2)}(0)&     \nu^{(2)}(1)&\cdots&\nu^{(2)}(s-1)
    \end{bmatrix}$ with
    \begin{equation}\label{eqn:coherence 2}
         \nu^{(2)}(t)=\max\limits_{1\le\ell\le n}
    \left\{\frac{1}{{\mathbf{p}^{(2)}_t(\ell)}}\left\|
   \delta_\ell^{\top}U_k
\begin{bmatrix}
    \Lambda_k^{t} & \bar{\Lambda}_k^{t}\\
\end{bmatrix}\right\|_2^2\right\}.
    \end{equation}
    \end{enumerate}
\end{definition}

Given the challenges in providing the closed-form expression for the exact value of $\nu^{(j)},$ we present the following estimate instead.

\begin{lemma}\label{nuestimate}
The graph weighted coherence $\nu^{(2)}$ is bounded by
\begin{equation*}
    \nu^{(2)}(t)\le \max_{1\leq i\leq n}\frac{\|U_k^\top\delta_i\|^2}{\mathbf{p}^{(2)}_t(i)} \max_{1\leq j\leq k}\left(\lambda_j^{2t}+(\sum_{l=1}^{t-1}\lambda_j^l)^2\right).
\end{equation*}
In particular, when $\mathbf{P}^{(j)}\ (j=1,2)$ are uniform distributions, $\nu^{(1)}\ge\nu^{(2)}(t)$ for all $t=0,1,\cdots,s-1.$
\end{lemma}

\subsection{Restricted isometry property}
Our goal is to recover the initial state  $\bx_0$  and the source term $\bw$  from a small number of space-time measurements. The Restricted Isometry Property (RIP) is a key concept in compressed sensing and signal recovery, ensuring that accurate recovery is possible when the measurements satisfy the RIP condition. In our work, we slightly relax this condition by requiring that, under certain conditions, the target signals can be accurately recovered with high probability.

After selecting the sampling node set \(\Omega_t = \{\omega_{1,t}, \cdots, \omega_{m_t,t}\}\), we define the weighting matrix as \(M_t = \frac{1}{\sqrt{m_t}}I \in \mathbb{R}^{m_t \times m_t}\), and the sampling distribution matrix as \(P_{\Omega_t} = \text{diag}(\mathbf{p}_t(\omega_{1,t}), \dots, \mathbf{p}_t(\omega_{m_t,t}))\). Additionally, we define the diagonal block matrices \(\mathcal{W}\), \(\mathcal{P}\), and \(\mathcal{S}\). A summary of these matrices, along with other related notations, is provided in \Cref{table:1}. We then present our main   \Cref{thm:smp_cmp} which shows that $\mathcal{W}\mathcal{P}\mathcal{S}\pi_{A,s}$ satisfies the RIP.

\begin{table}[ht]\label{def}
    \centering
    \renewcommand{\arraystretch}{1.6}
    \small
    \begin{tabular}{|c|c|}
        \hline
        \textbf{Notation} & \textbf{Definition} \\ 
        \hline
        $\Omega_t$ & $\Omega_t=\{\omega_{1,t},\cdots, \omega_{m_t,t}\}$ \\ 
        \hline
        $M_t$ &  $M_t=\frac{1}{\sqrt{m_t}}I \in \mathbb{R}^{m_t \times m_t}$ \\ 
        \hline
        $P_{\Omega_t}$ &  $P_{\Omega_t}=\text{diag}(\mathbf{p}_t(\omega_{1,t}),\dots,\mathbf{p}_t(\omega_{m_t,t}))\in \mathbb{R}^{m_t \times m_t}$ \\ 
        \hline
        $\mathcal{W}$ & 
        $\mathcal{W}=\begin{bmatrix}
            M_0 &  &  &  \\
            & M_1 &  &  \\
            &  & \ddots &  \\
            &  &  & M_{s-1}
        \end{bmatrix} \in \mathbb{R}^{\sum\limits_{t=0}^{s-1}m_t \times \sum\limits_{t=0}^{s-1}m_t}$ \\ 
        \hline
        $\mathcal{P}$ & 
        $\mathcal{P}=\begin{bmatrix}
            P_{\Omega_0}^{-\frac{1}{2}} &  &  &  \\
            & P_{\Omega_1}^{-\frac{1}{2}} &  &  \\
            &  & \ddots &  \\
            &  &  & P_{\Omega_{s-1}}^{-\frac{1}{2}}
        \end{bmatrix} \in \mathbb{R}^{\sum\limits_{t=0}^{s-1}m_t \times \sum\limits_{t=0}^{s-1}m_t}$ \\ 
        \hline
        $\mathcal{S}$ & 
        $\mathcal{S}=\begin{bmatrix}
            S_0 &  &  &  \\
            & S_1 &  &  \\
            &  & \ddots &  \\
            &  &  & S_{s-1}
        \end{bmatrix} \in \mathbb{R}^{\sum\limits_{t=0}^{s-1}m_t \times Ns}$ \\ 
        \hline
        $\nu^{(1)}$ & $\nu^{(1)}\in\mathbb{R}$ in  \Cref{coherence}\\
        \hline
        $\nu^{(2)}$ &
        $\nu^{(2)}=\begin{bmatrix}
        \nu^{(2)}(0)&\cdots&\nu^{(2)}(s-1)
    \end{bmatrix}$ in  \Cref{coherence}\\
    \hline
    \end{tabular}
    \caption{Definitions of matrices and notations}
    \label{table:1}
\end{table}

\begin{theorem}\label{thm:smp_cmp}
    Let $\bw_{\text{aug}}=\begin{bmatrix}
   \bx_0\\
   \bw
\end{bmatrix}$ and  $\mathcal{W}, \mathcal{P}, \mathcal{S}$ be defined in \Cref{table:1}   with the sampling distribution $\mathbf{P} = [\mathbf{p}_0,\mathbf{p}_1, \cdots, \mathbf{p}_{s-1}].$ For any $\delta, \epsilon\in(0,1),$ with probability at least $1-\epsilon,$
\begin{equation}\label{RIP}
(1-\delta)c_{A,k,s}\left\|
\bw_{\text{aug}}\right\|^2
\le
\left\|\mathcal{W}\mathcal{P}\mathcal{S}\pi_{A,s}(\bw_{\text{aug}})\right\|^2\le(1+\delta)C_{A,k,s}\left\|\bw_{\text{aug}}\right\|^2
\end{equation}
   for all $\bw_{\text{aug}}\in\text{span}(\widetilde{U}_k)$ provided that  
   \[
m\ge \frac{3\nu^{(1)}}{c_{A,k,s}\delta^2}\log\frac{4k}{\epsilon}
\] in regime 1, and 
\[
m_t\ge \frac{3\nu^{(2)}(t)}{c_{A,k,s}\delta^2}\log\frac{4k}{\epsilon}
\]
for $t=0,\cdots,s-1$ in regime 2. 

\end{theorem}

\section{Reconstruction Algorithms}\label{algorithm}

In this section, we aim to design a robust procedure to recover the initial signal $\bx_0$ and the source term $\bw$ using $\sum\limits_{t=0}^{s-1}m_t$ space-time samples. More specifically, the space-time samples $\bz\in\mathbb{R}^{\sum\limits_{t=0}^{s-1}m_t}$ are given by
\begin{equation}
    \bz=\mathcal{S}\pi_{A,s}({\bw}_{\text{aug}})+\be
\end{equation}
where $\bw_{\text{aug}}=\begin{bmatrix}
    \bx_0\\
    \bw
\end{bmatrix}$ and $\be\in\R^{\sum\limits_{t=0}^{s-1}m_t}$ is the observational noise. When a basis for $\text{span}(U_k)$ is known, we use the standard least squares method to estimate $\bw_{\text{aug}}$ from $\bz$ by solving
\begin{equation}\label{recons1}    \bw_{\text{aug}}^*=\arg\min_{\widetilde{\bw}_{\text{aug}}\in\text{span}(\Tilde{U}_k)}\|\mathcal{W}\mathcal{P}(\mathcal{S}\pi_{A,s}(\widetilde{\bw}_{\text{aug}})-\bz)\|_2,
\end{equation} where
$\widetilde{U}_k=\begin{bmatrix}
    U_k &\textbf{0}\\
    \textbf{0} & U_k
\end{bmatrix}$. One can solve \eqref{recons1} via solving
\begin{equation}\label{recons1-1}    \bv^*=\arg\min_{\overline{\bv}\in\R^{2k}}\|\mathcal{W}\mathcal{P}(\mathcal{S}\pi_{A,s}(\widetilde{U}_k\overline{\bv})-\bz)\|_2
\end{equation}
and setting $\bw_{\text{aug}}^*=\widetilde{U}_k\bv^*.$

\begin{theorem}\label{alg1}
The operators \(\mathcal{W}, \mathcal{P},\) and \(\mathcal{S}\) are defined in Table \ref{table:1}, while \(c_{A,k,s}\) and \(C_{A,k,s}\) represent the bounds of the operator \(\pi_{A,s}\). Let $(\epsilon, \delta)\in(0,1),$ and suppose that $
m\ge \frac{3\nu^{(1)}}{c_{A,k,s}\delta^2}\log\frac{4k}{\epsilon}
$ in regime 1 and 
$m_t\ge \frac{3\nu^{(2)}(t)}{c_{A,k,s}\delta^2}\log\frac{4k}{\epsilon}$ for $t=0,\cdots,s-1$ in regime 2. With probability at least $1-\epsilon,$ the following holds for all $\bw_{\text{aug}}\in \text{span}(\widetilde{U}_k)$ and all $\be\in\R^{\sum\limits_{t=0}^{s-1}m_t}.$ Then
\begin{enumerate}
    \item Let $\bw_{\text{aug}}^*$ be the solution of \eqref{recons1} with space-time samples $\bz=\mathcal{S}\pi_{A,s}(\bw_{\text{aug}})+\be,$ then \begin{equation}\label{noiseub}
    \left\|
\bw_{\text{aug}}^*-\bw_{\text{aug}}\right\|_2\le \frac{2}{\sqrt{(1-\delta)c_{A,k,s}}}\|\mathcal{W}\mathcal{P}\be\|_2.
\end{equation}.\\
    \item There exists some particular vector $\widetilde{\be}\in\R^{\sum\limits_{t=0}^{s-1}m_t}$ such that the solution of \eqref{recons1} with space-time samples $\bz=\mathcal{S}\pi_{A,s}(\bw_{\text{aug}})+\widetilde{\be}$ satisfies 
    \begin{equation}\label{noiselb}
    \|\bw_{\text{aug}}^*-\bw_{\text{aug}}\|_2\ge
    \frac{1}{\sqrt{(1+\delta)C_{A,k,s}}}\|\mathcal{W}\mathcal{P}\widetilde{\be}\|_2.
\end{equation}
\end{enumerate}
\end{theorem}
The estimate \eqref{noiseub} implies that when \( \be = \textbf{0} \), we have \( \bw_{\text{aug}}^* = \bw_{\text{aug}} \) with high probability. In the presence of noise, \eqref{noiseub} establishes that \( \|\bw_{\text{aug}}^* - \bw_{\text{aug}}\|_2 \) grows linearly with \( \|\mathcal{W}\mathcal{P}\be\|_2 \).

Notice that solving \eqref{recons1} requires prior knowledge of $U_k,$ which can be computationally expensive for large scale graphs. When $U_k$ is unavailable, we propose estimating $\bw_{\text{aug}}$ by solving the following regularized least squares problem.
\begin{equation}\label{recons2}   \bw^*_{\text{aug}} =  \arg\min_{\overline{\bv}\in\R^{2n}}\|\mathcal{W}\mathcal{P}(\mathcal{S}\pi_{A,s}(\overline{\bv})-\bz)\|^2_2+\gamma\overline{\bv}^{\top}g(\widetilde{L})\overline{\bv}
\end{equation}
where $\gamma>0,$ $\widetilde{L}=\begin{bmatrix}
    L & \textbf{0}\\
    \textbf{0} & L
\end{bmatrix},$ and $g:\mathbb{R}\rightarrow\mathbb{R}$ is a nonnegative and non-decreasing polynomial function. The term $\gamma\overline{\bv}^{\top}g(\widetilde{L})\overline{\bv}$ is used to penalize signals that have significant
high frequency components. \eqref{recons2} can be solved by solving the following equation
\[
(\pi_{A,s}^\top S^\top\mathcal{P}^\top\mathcal{W}^\top\mathcal{W}\mathcal{P}S\pi_{A,s}+\gamma g(\widetilde{L}))\overline{\bv}=\pi_{A,s}^\top S^\top\mathcal{P}^\top\mathcal{W}^\top\mathcal{W}\mathcal{P}\bz.
\]
The next theorem provides bounds for the error between the original signal $\bw_{\text{aug}}$ and the solution of the optimization problem \eqref{recons2}.

\begin{theorem}\label{thm:alg2}
The operators \(\mathcal{W}, \mathcal{P},\) and \(\mathcal{S}\) are defined in Table \ref{table:1}, while \(c_{A,k,s}\) and \(C_{A,k,s}\) represent the bounds of the operator \(\pi_{A,s}\) defined in \eqref{boundpi}. $R$ is a constant such that $\|\mathcal{W}\mathcal{P}\mathcal{S}\|_2\le R.$ Let $(\epsilon, \delta)\in(0,1),$ and suppose that $
m\ge \frac{3\nu^{(1)}}{c_{A,k,s}\delta^2}\log\frac{4k}{\epsilon}
$ in regime 1 and 
$m_t\ge \frac{3\nu^{(2)}(t)}{c_{A,k,s}\delta^2}\log\frac{4k}{\epsilon}$ for $t=0,\cdots,s-1$ in regime 2. With probability at least $1-\epsilon,$ the following holds for all $\bw_{\text{aug}}\in \text{span}(\widetilde{U}_k)$ and all $\be\in\R^{\sum\limits_{t=0}^{s-1}m_t},$ all $\gamma>0,$ and all non-negative, non-decreasing polynomial functions $g$ such that $g(\theta_{k+1})>0.$ 

Let $\bw_{\text{aug}}^*$ be the solution of \eqref{recons2} with space-time samples $\bz=\mathcal{S}\pi_{A,s}(\bw_{\text{aug}})+\be,$ then 
\begin{equation}\label{alphabd}
        \|\alpha^*-\bw_{\text{aug}}\|_2\le\left(\frac{2+\frac{RC_{A,n,s}}{\sqrt{\gamma g(\theta_{k+1})}}}{\sqrt{(1-\delta)c_{A,k,s}}}\right)\|\mathcal{W}\mathcal{P}\be\|_2
    +\left(\frac{RC_{A,n,s}\sqrt{\frac{g(\theta_{k})}{g(\theta_{k+1})}}+\sqrt{\gamma g(\theta_k)}}{\sqrt{(1-\delta)c_{A,k,s}}}\right)\|\bw_{\text{aug}}\|_2.
\end{equation}
and
\begin{equation}\label{betabd}
\|\beta^*\|_2\le\frac{\|\mathcal{W}\mathcal{P}\be\|_2}{\sqrt{\gamma g(\theta_{k+1})}}
+\sqrt{\frac{g(\theta_k)}{g(\theta_{k+1})}}\|\bw_{\text{aug}}\|_2
\end{equation}
where $\alpha^*=\widetilde{U}_k\widetilde{U}_k^\top\bw_{\text{aug}}^*$ and $\beta^*=(I_{2n}-\widetilde{U}_k\widetilde{U}_k^\top)\bw_{\text{aug}}^*.$
\end{theorem}

\begin{remark}
According to the bounds in \eqref{alphabd} and \eqref{betabd}, the error \(\|\bw_{\text{aug}}^* - \bw_{\text{aug}}\|\) can be estimated as
\begin{equation}
    \|\bw_{\text{aug}}^* - \bw_{\text{aug}}\|^2 = \|\alpha^* - \bw_{\text{aug}}\|^2 + \|\beta^*\|^2 \leq \left( \|\alpha^* - \bw_{\text{aug}}\| + \|\beta^*\| \right)^2.
\end{equation}
In particular, in the absence of noise, we have the following bound:
\begin{equation}
    \|\bw_{\text{aug}}^* - \bw_{\text{aug}}\| \le \left( \frac{R C_{A,n,s} \sqrt{ \frac{g(\theta_k)}{g(\theta_{k+1})} } + \sqrt{ \gamma g(\theta_k) } }{ \sqrt{(1 - \delta) c_{A,k,s}} } + \sqrt{ \frac{g(\theta_k)}{g(\theta_{k+1})} } \right) \|\bw_{\text{aug}}\|_2.
\end{equation}
This expression indicates that if \(g\) is chosen such that \(g(\theta_k) = 0\), then perfect reconstruction is achievable. In the more general case where \(g(\theta_k) > 0\), good reconstruction performance can still be expected provided that the ratio \(\frac{g(\theta_k)}{g(\theta_{k+1})}\) is kept sufficiently small and the regularization parameter \(\gamma\) is close to zero.
\end{remark}

\section{Experimental results}\label{experiments}
This section introduces the datasets used in this work and details the experimental results. The experiments involve four datasets: two synthetic datasets and two real-world datasets. For the synthetic data, we use the same set of graphs as in \cite{huang2024robust}, comprising two distinct graph types available in the GSP toolbox \cite{perraudin2014gspbox}:
(a) the Minnesota road graph with \( n = 2642 \) nodes; and (b) the Stanford bunny graph with \( n = 2503 \) nodes. The real-world datasets include (c) the sea surface temperature recorded monthly and provided by the NOAA Physical Sciences Laboratory (PSL), and (d) the global COVID-19 dataset
provided by the Johns Hopkins University \cite{dondugar2005}. 


\paragraph{\textbf{A. Datasets}}
\begin{enumerate}
    \item Synthetic datasets: We use the synthetic datasets, which was updated by Huang et al. \cite{huang2024robust}, based on the original datasets provided in the GSP toolbox \cite{perraudin2014gspbox}. Motivated by the sea surface temperature data, we model continuous heat diffusion on all synthetic datasets using the following formulation: 
    \[\bx_{t+1}=A\bx_t+\bw\] where 
    \( A := e^{-\alpha L} \), and \( L \) represents the normalized Laplacian matrix. The parameter $\alpha$ varies across different cases to capture the dynamics of heat flow on the graph structure. 
    In our experiment, we set \(\alpha = 30\). For all synthetic graphs, we choose \(s = 10\), which represents the total sampling time.\\
    \item Sea Surface Temperature: The graph \( G \) is constructed using the \( k \)-nearest neighbors (k-NN) algorithm with $k=10$, based on the coordinate locations of the nodes in each dataset, as described in \cite{9730033}. That is, let \( R \in \mathbb{R}^{n \times 2} \) be the matrix of coordinates for all nodes, such that \( R = [r_1,r_2, \dots, r_n]^\top \), where \( r_i \in \R^2 \) is the vector containing the latitude and longitude of node \( i \). The weight of each edge \( (i,j) \) is given by the Gaussian kernel:
\[
W(i,j) = \exp\left( - \frac{\| r_i - r_j \|_2^2}{\sigma^2} \right),
\]
where \( \| r_i - r_j \|_2^2 \) is the squared Euclidean distance between nodes \( i \) and \( j \), and \( \sigma \) is the standard deviation, defined as
\[
\sigma = \frac{1}{|\mathcal{E}|} \sum_{(i,j) \in \mathcal{E}} \| r_i - r_j \|_2.
\]
Here, \( \mathcal{E} \) represents the set of edges in the graph. Furthermore, \( A := e^{-\alpha L} \), where the combinatorial Laplacian is used. In this dataset, we use a subset of 100 points on the Pacific Ocean
within a time frame of 110 months. The first 10 months are used for training, while the following 100 months are used for experiments.
\item Global COVID-19: This dataset contains the cumulative number of daily COVID-19 cases for 265 locations worldwide. The weighted adjacency matrix is defined using the same method as described for the sea surface temperature dataset, and $A$ is also defined similarly. We use the data from January 22, 2020, covering a total of 110 days. Likewise, the first 10 days are used for training, while the following 100 days are used for experiments.
\end{enumerate}

\paragraph{\textbf{B. Evaluation Metrics}} 
We use the following evaluation metrics: 
Mean Absolute Error (MAE), Mean Absolute Percentage Error (MAPE), and Relative Error (RE). These metrics are defined as follows:


\[\text{MAE} = \frac{1}{n_x} \sum_{i=1}^{n_x} |\mathbf{x}^*(i) - \mathbf{x}(i)|,\quad\text{MAPE} = \frac{1}{n_x} \sum_{i=1}^{n_x} \left| \frac{\mathbf{x}^*(i) - \mathbf{x}(i)}{\mathbf{x}(i)} \right|,\quad\text{RE} = \frac{\|\mathbf{x}^*-\mathbf{x}\|_2}{\|\mathbf{x}\|_2}
\]
where \( \mathbf{x}^*(i) \) is the \( i \)-th element of the recovered signal, \( \mathbf{x}(i) \) is the \( i \)-th element of the ground truth signal, and \( n_x \) is the length of the signal.

\subsection{Experiments with Synthetic Data} In the experiments with synthetic data, we first generate the synthetic vector \(\bw_{\text{aug}} = \begin{bmatrix} \bx_0 \\ \bw \end{bmatrix} \in \widetilde{U}_k\) and then perform separate reconstructions for both regime 1 and regime 2 using \eqref{recons1} and \eqref{recons2}. We evaluate the reconstruction accuracy using relative error and investigate the impact of noise and the parameter \(\gamma\). While our theoretical framework is applicable to any probability distribution for sampling graph nodes, we specifically use uniform probability in our experiments.


The total number of samples is given by \( M,\) and all space-time samples are drawn with replacement. In  regime 1, the sampling locations are fixed over time, and each time step uses \( m_t = M/s \) samples for all \( t \in [s] \). In regime 2, the sampling locations are selected independently at each time step, while maintaining \( m_t = M/s \) for all \( t \in [s].\) Recall that \( s = 10 \) in these experiments.

\begin{figure}[ht]
    \centering
    \begin{minipage}{0.45\linewidth}
        \centering
    \includegraphics[width=\linewidth]{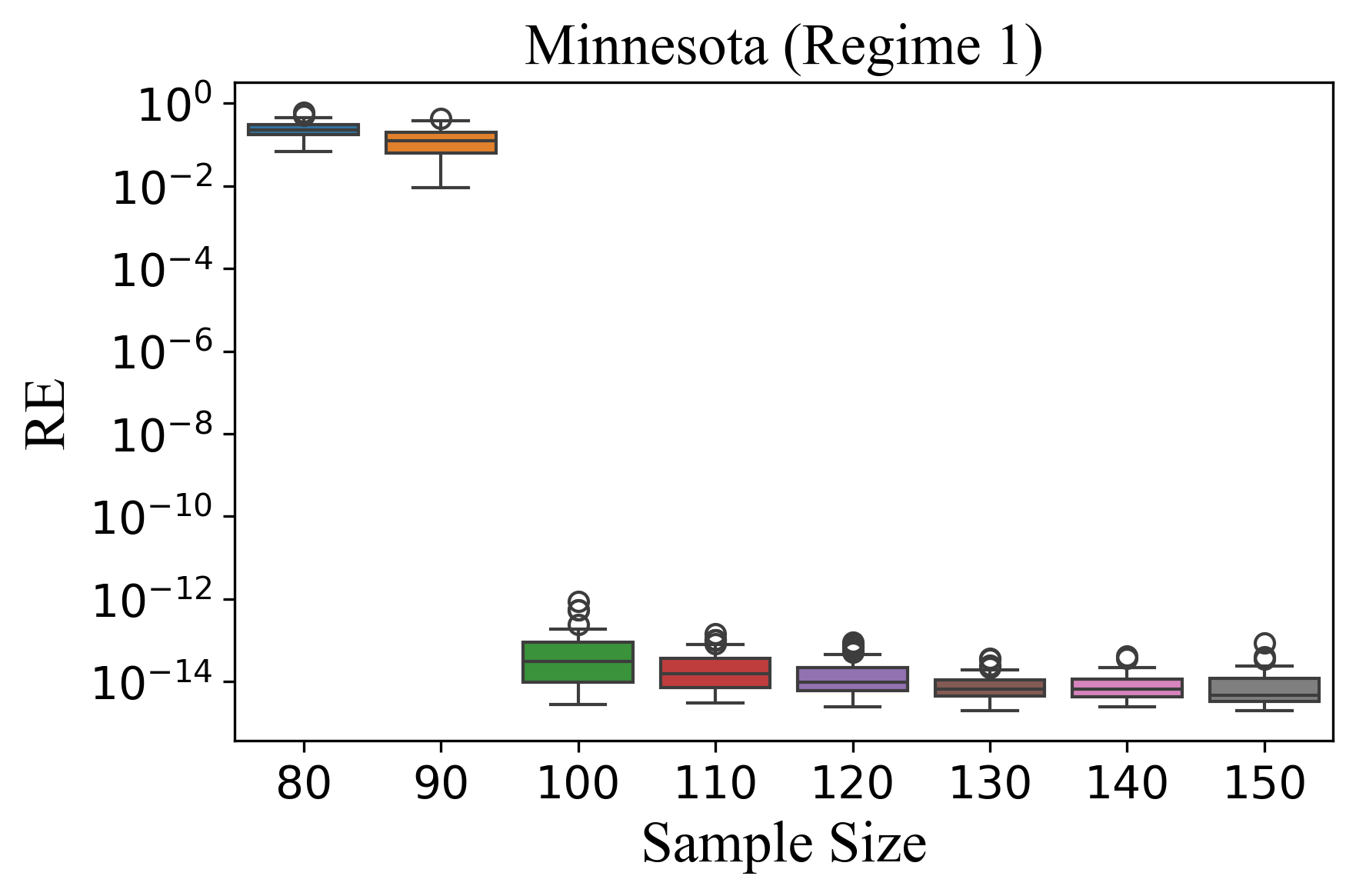}
    \end{minipage}
    \begin{minipage}{0.45\linewidth}
        \centering
    \includegraphics[width=\linewidth]{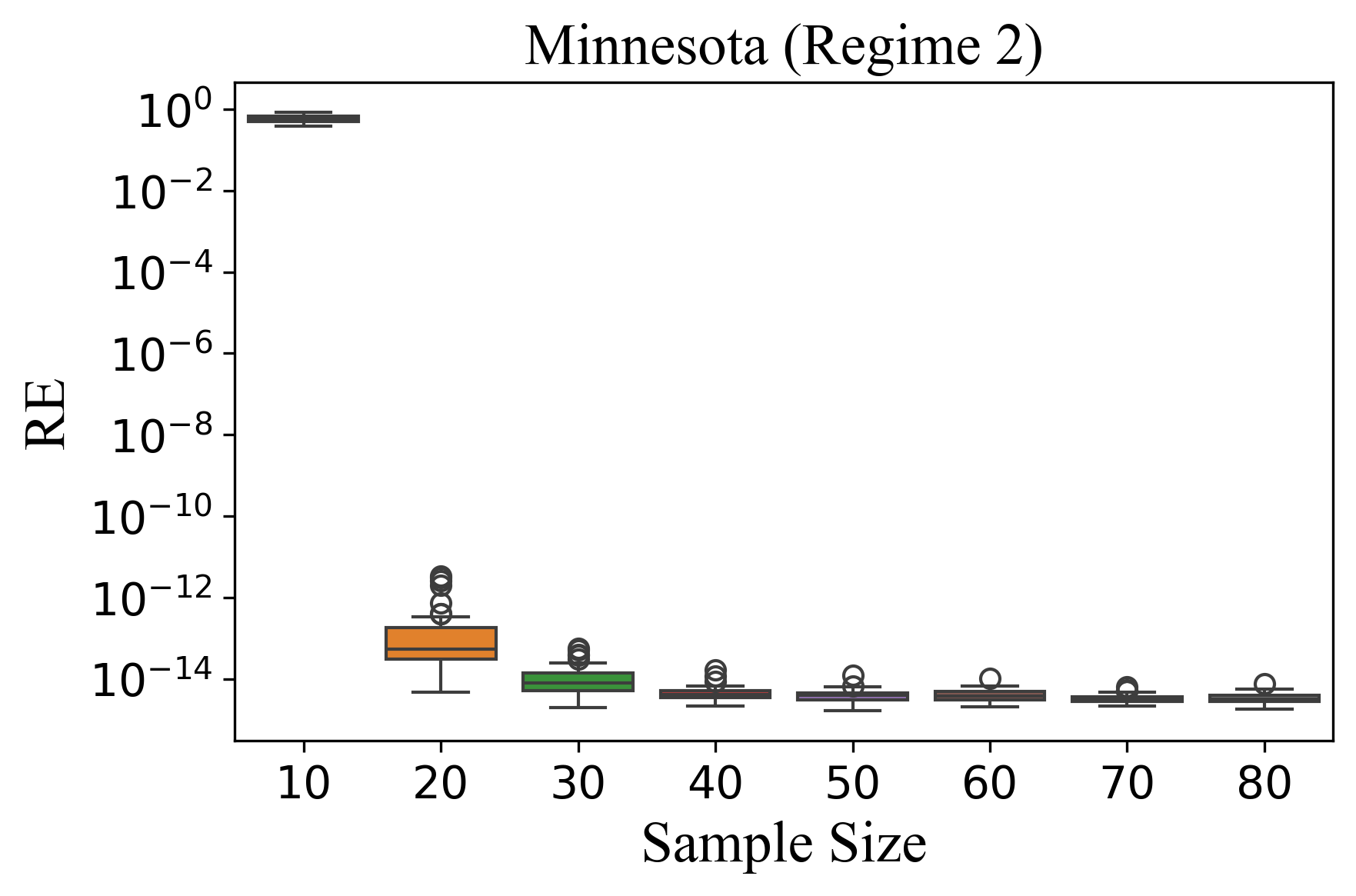}
    \end{minipage}
    
    \begin{minipage}{0.45\linewidth}
        \centering
    \includegraphics[width=\linewidth]{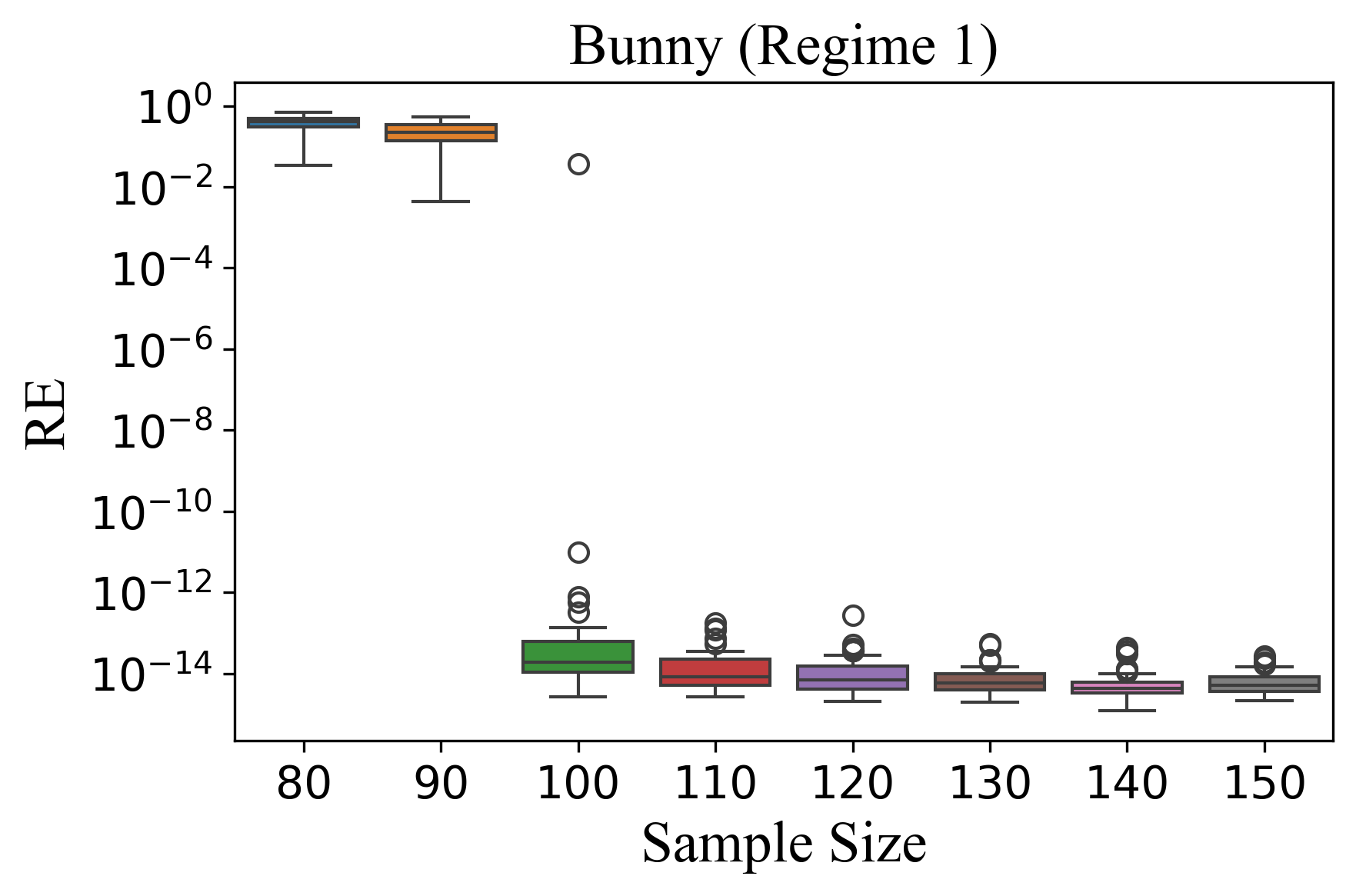}
    \end{minipage}
    \begin{minipage}{0.45\linewidth}
        \centering
    \includegraphics[width=\linewidth]{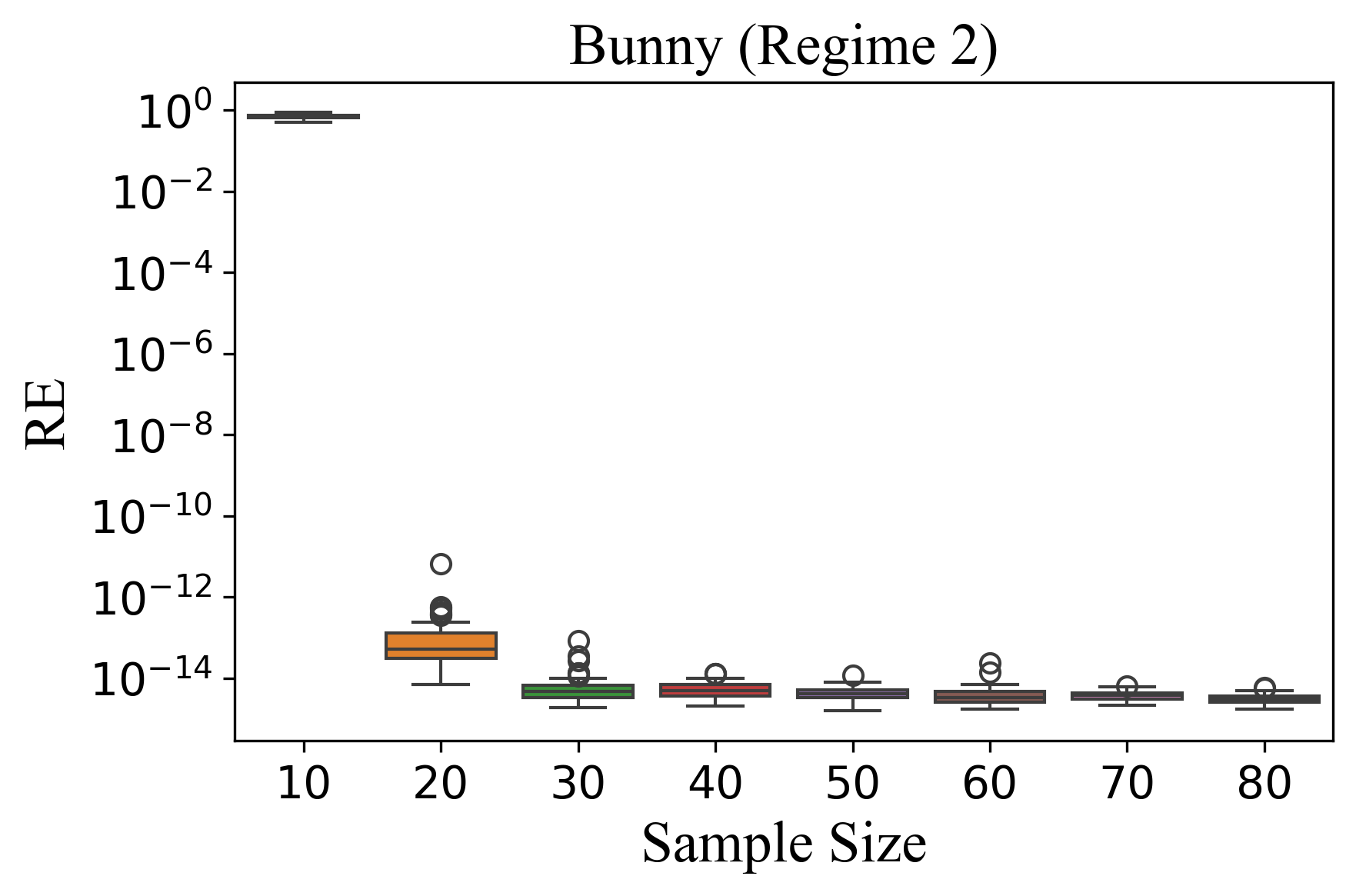}
    \end{minipage}

    \caption{Relative error vs. sample size \( M \) without noise. The first row corresponds to the Minnesota graph, while the second row corresponds to the Bunny graph. The first column shows results for regime 1, and the second column shows results for regime 2.}
    \label{fig1}
\end{figure}

\begin{figure}[ht]
    \centering
    \begin{minipage}{0.45\linewidth}
        \centering
        \includegraphics[width=\linewidth]{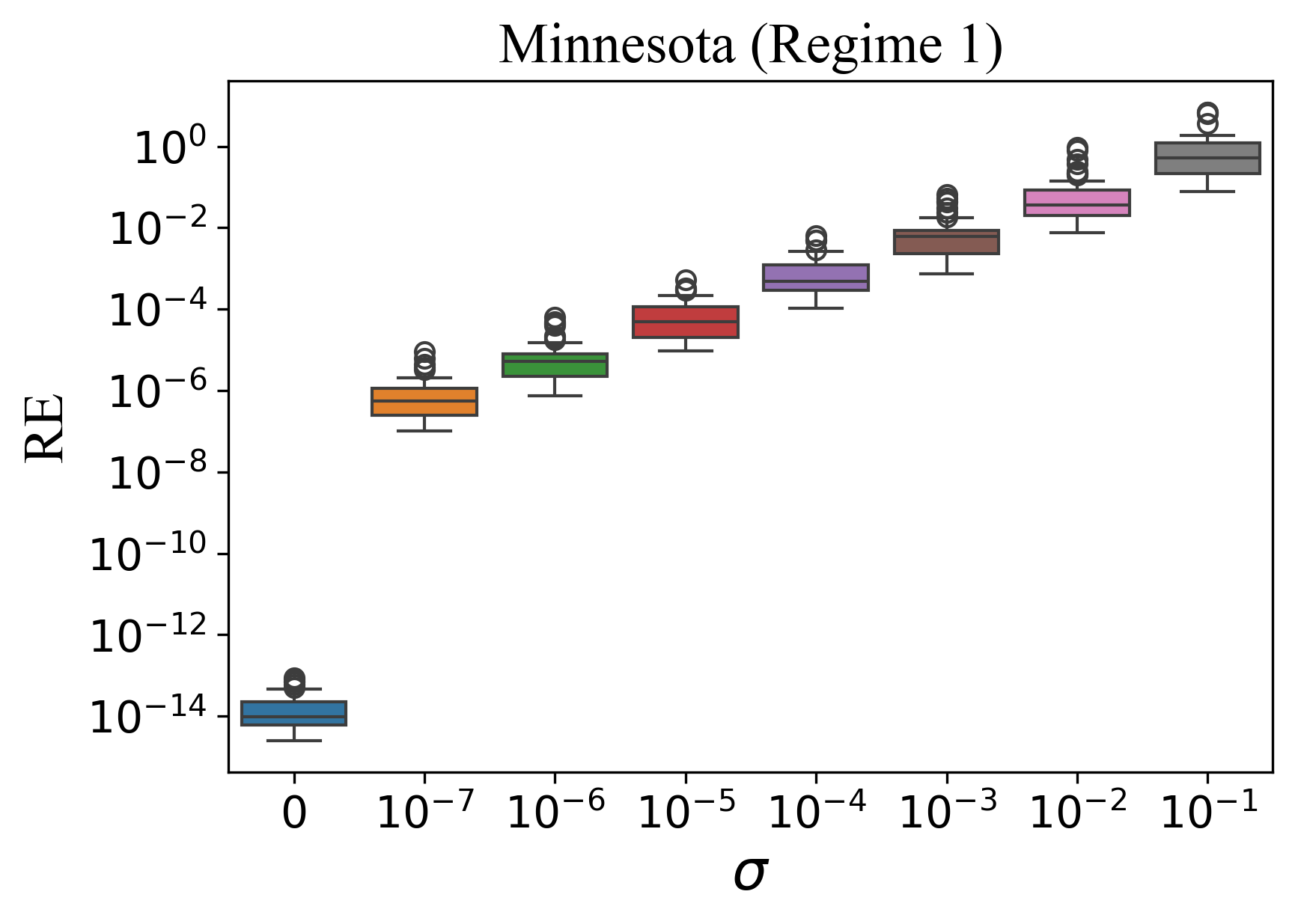}
    \end{minipage}
    \begin{minipage}{0.45\linewidth}
        \centering
        \includegraphics[width=\linewidth]{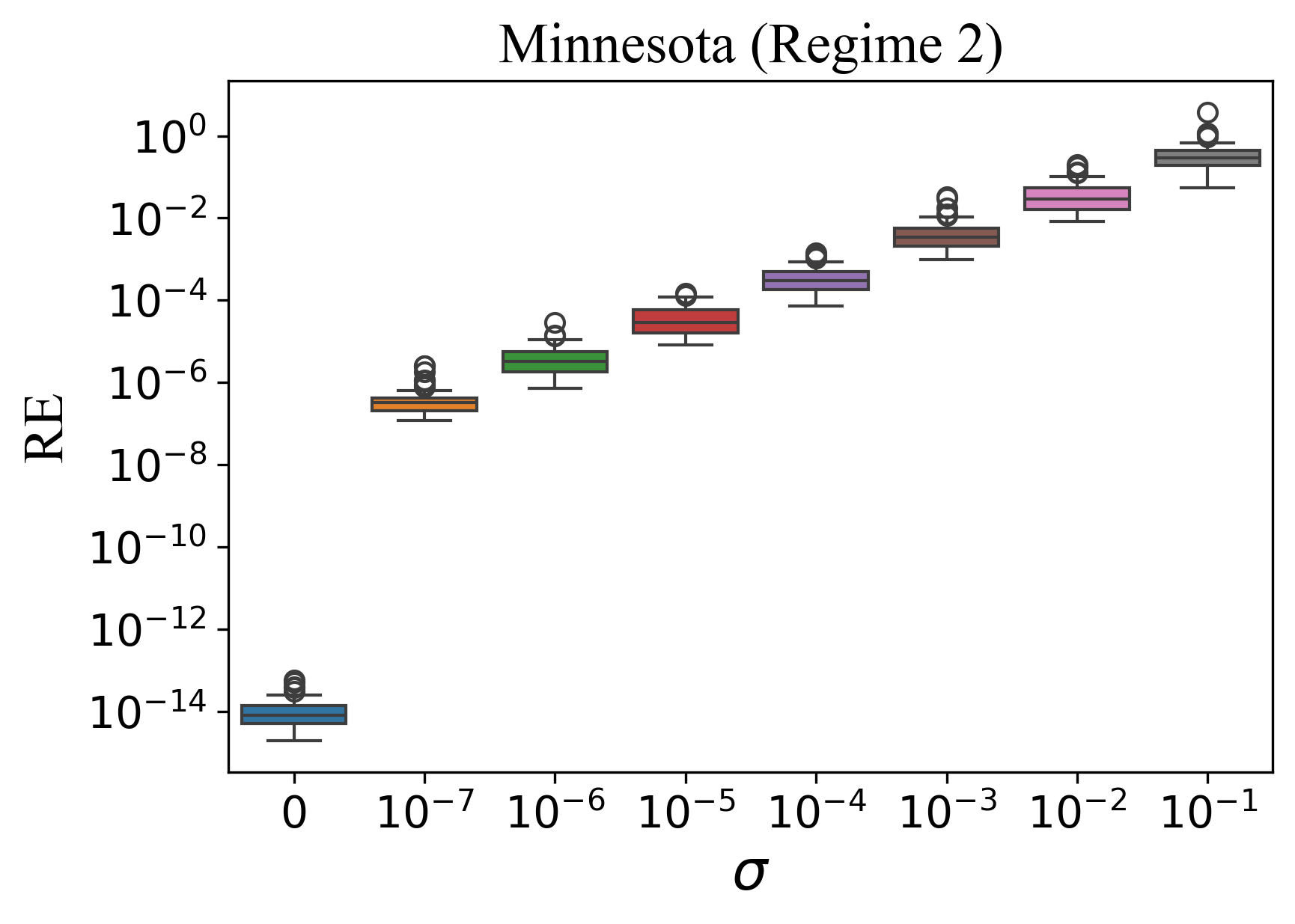}
    \end{minipage}
    
    
    \begin{minipage}{0.45\linewidth}
        \centering
        \includegraphics[width=\linewidth]{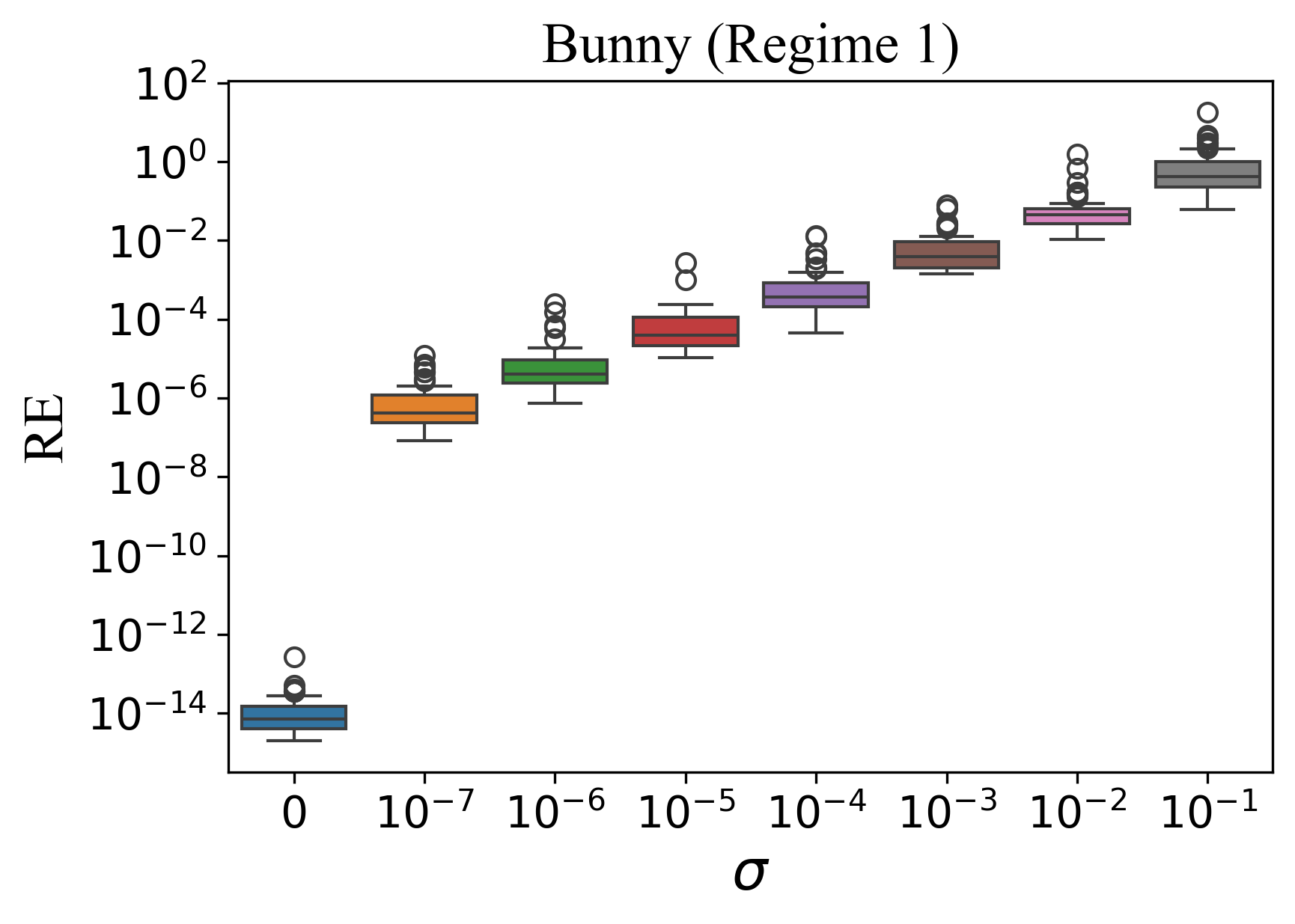}
    \end{minipage}
    \begin{minipage}{0.45\linewidth}
        \centering
        \includegraphics[width=\linewidth]{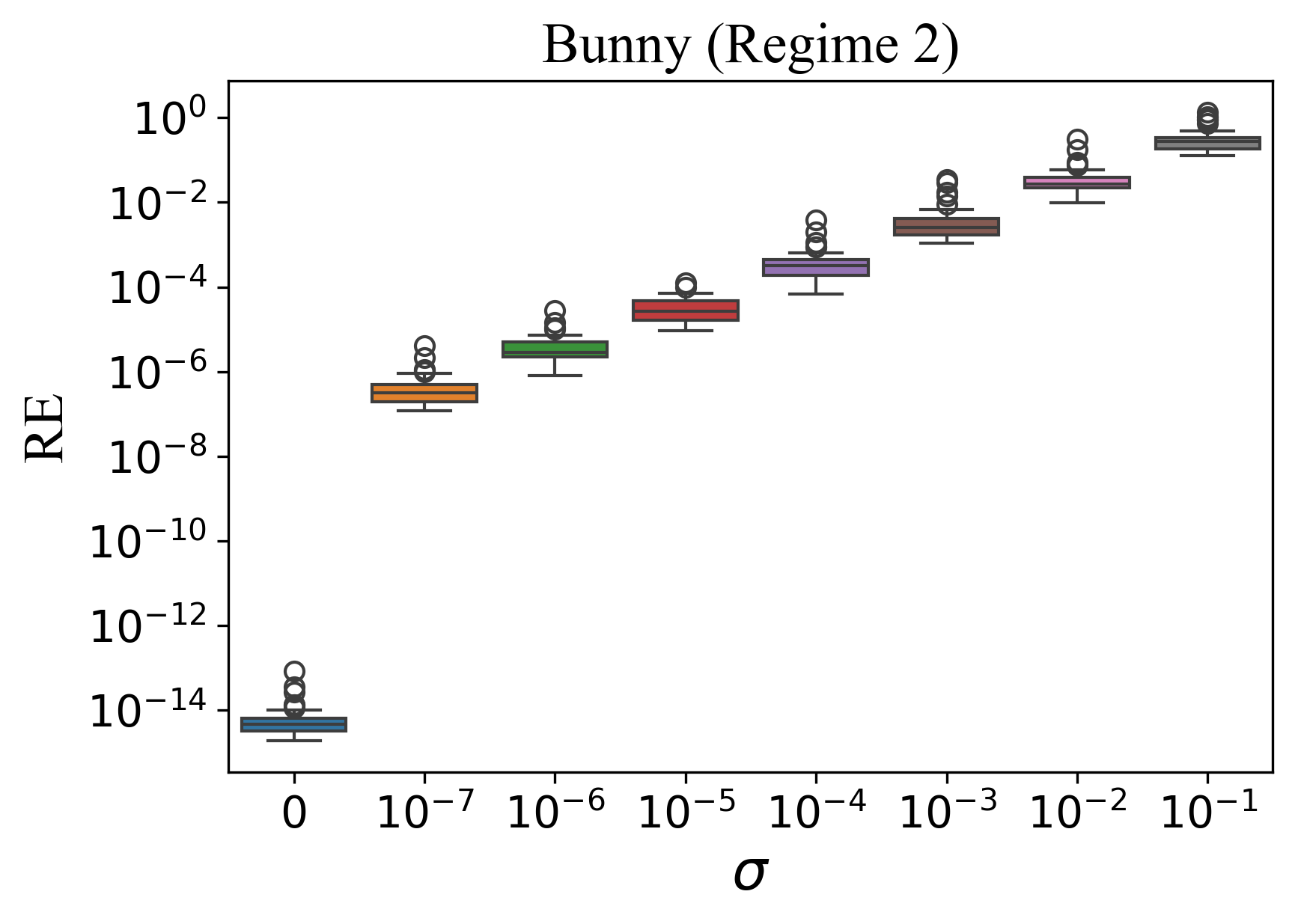}
    \end{minipage}

    \caption{Relative error vs. noise level. The first row corresponds to the Minnesota graph, while the second row corresponds to the Bunny graph. The first column shows results for regime 1, and the second column shows results for regime 2.}
    \label{fig2}
\end{figure}

\subsubsection{Reconstruction with $U_k$}\label{nonreg} When \( U_k \) is assumed to be known, we use \eqref{recons1} to recover the vector \( \bw_{\text{aug}} \). The reconstruction results without noise for both regime 1 and regime 2 are presented in Fig.~\ref{fig1}. For each fixed \( M \), we conducted 50 experiments and generated a boxplot to visualize the results. The experiment demonstrates that 100 total samples are sufficient to recover \( \bw_{\text{aug}}  \) in regime 1, while 20 samples are sufficient in regime 2 for both the Minnesota and Bunny graphs. It demonstrates that the flexibility in selecting sampling nodes in regime 2 leads to a more efficient use of the total samples, achieving the same performance with fewer samples compared to regime 1.

To evaluate the impact of noise on reconstruction, we first fix the sample size \( M \). We generate a random noise vector whose entries follow a normal distribution with standard deviation \( \sigma \). Experiments are conducted with \( M = 110 \) in regime 1 and \( M = 30 \) in regime 2. We vary \( \sigma \) over the set \(\{0, 10^{-7}, 10^{-6}, 10^{-5}, 10^{-4}, 10^{-3}, 10^{-2}, 10^{-1} \}\), repeat the experiment 50 times for each setting, and visualize the results using a boxplot. The results are shown in \Cref{fig2}. As predicted by our theoretical analysis, the reconstruction error increases approximately linearly with respect to \( \sigma \).

\begin{figure}[ht]
    \centering
    \begin{minipage}{0.32\linewidth}
        \centering
      \includegraphics[width=\linewidth]{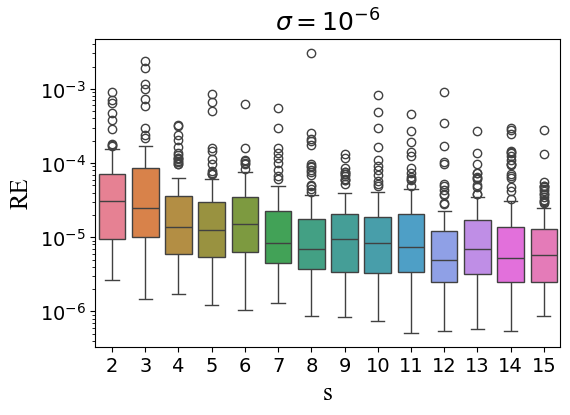}
    \end{minipage}
    \begin{minipage}{0.32\linewidth}
        \centering
   \includegraphics[width=\linewidth]{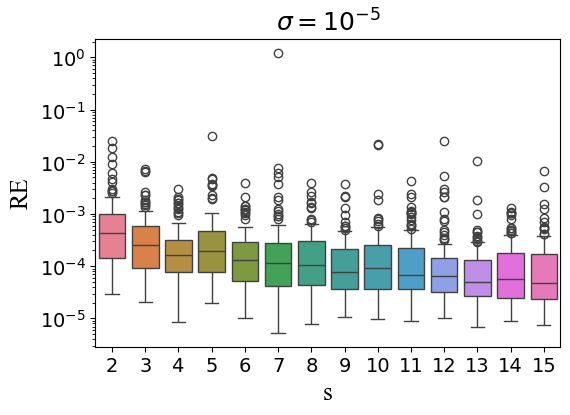}
    \end{minipage}
    \begin{minipage}{0.32\linewidth}
        \centering
      \includegraphics[width=\linewidth]{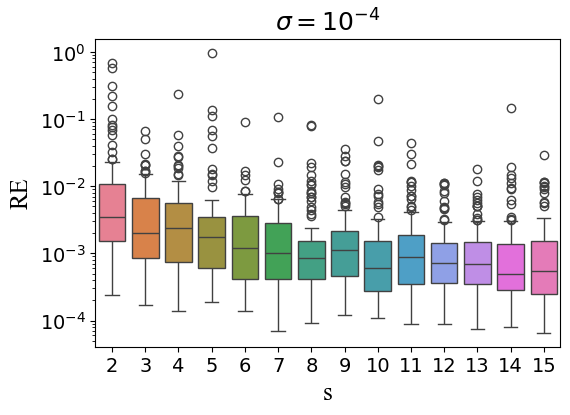}
    \end{minipage}
    
    
    \begin{minipage}{0.32\linewidth}
        \centering
    \includegraphics[width=\linewidth]{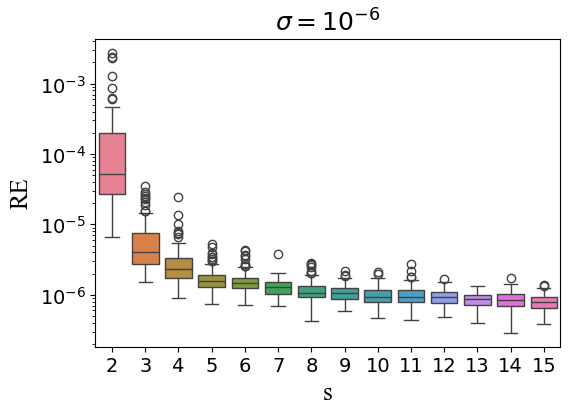}
    \end{minipage}
    \begin{minipage}{0.32\linewidth}
        \centering
        \includegraphics[width=\linewidth]{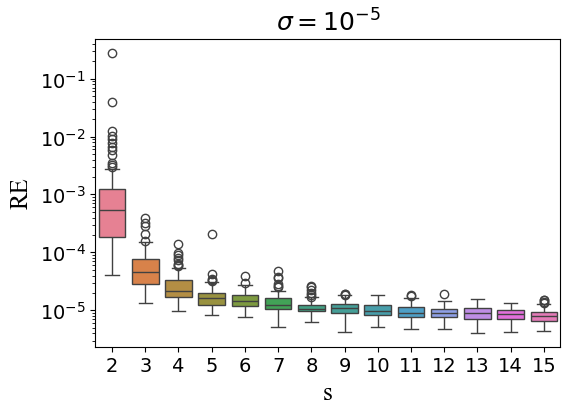}
    \end{minipage}
    \begin{minipage}{0.32\linewidth}
        \centering
        \includegraphics[width=\linewidth]{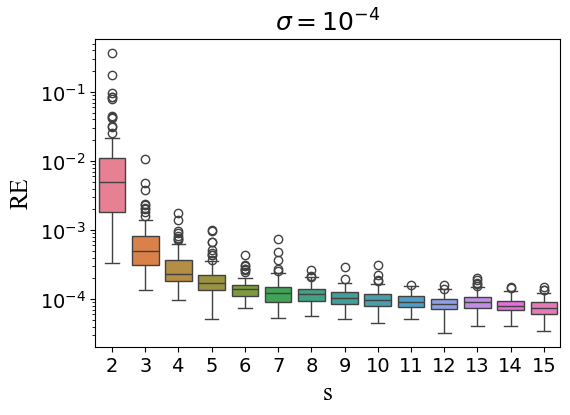}
    \end{minipage}
    
    \caption{Relative error vs. sample time $s$ at different noise levels for the Minnesota graph. The first row corresponds to regime 1, and the second row corresponds to regime 2.}
    \label{fig_s_minnesota}
\end{figure}

\begin{figure}[ht]
    \centering
    \begin{minipage}{0.32\linewidth}
        \centering
    \includegraphics[width=\linewidth]{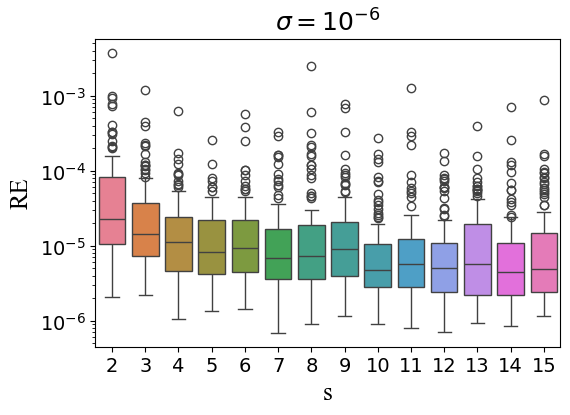}
    \end{minipage}
    \begin{minipage}{0.32\linewidth}
        \centering
    \includegraphics[width=\linewidth]{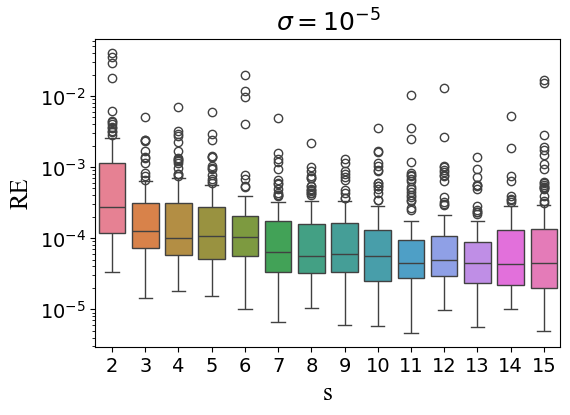}
    \end{minipage}
    \begin{minipage}{0.32\linewidth}
        \centering
    \includegraphics[width=\linewidth]{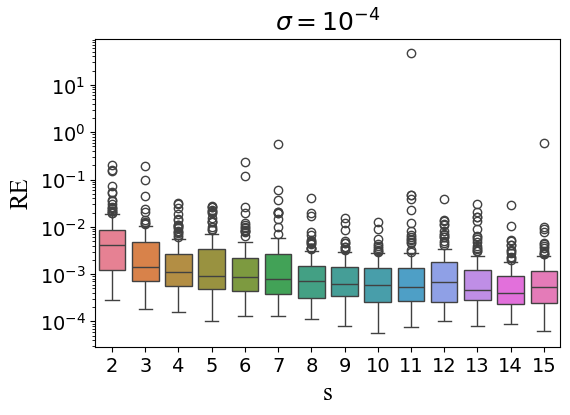}
    \end{minipage}
    
    
    \begin{minipage}{0.32\linewidth}
        \centering
    \includegraphics[width=\linewidth]{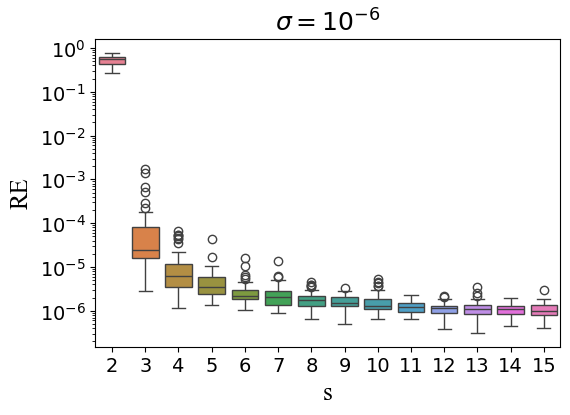}
    \end{minipage}
    \begin{minipage}{0.32\linewidth}
        \centering
    \includegraphics[width=\linewidth]{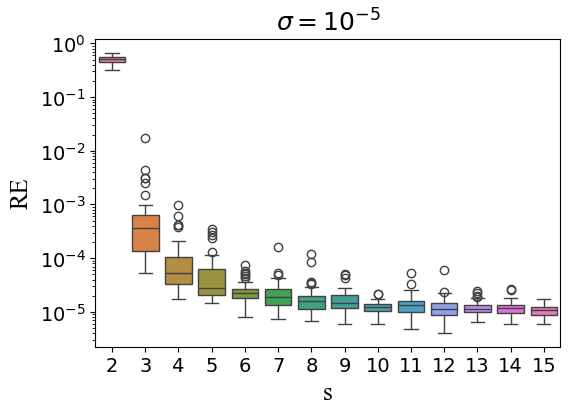}
    \end{minipage}
    \begin{minipage}{0.32\linewidth}
        \centering
    \includegraphics[width=\linewidth]{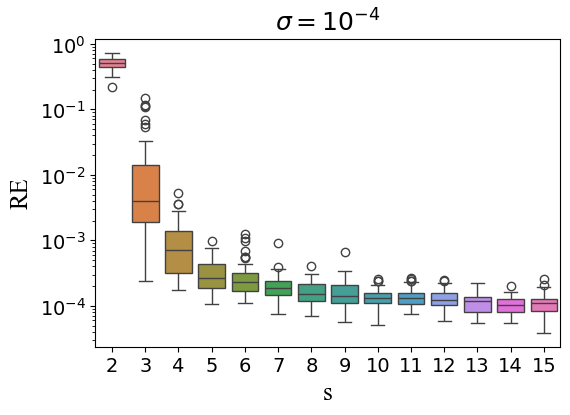}
    \end{minipage}
    \caption{Relative error vs. sample time $s$ at different noise levels for the Bunny graph. The first row corresponds to regime 1, and the second row corresponds to regime 2.}
    \label{fig_s_bunny}
\end{figure}

In the final experiment of this subsection, we investigate the effect of the total number of time steps $s$ on reconstruction accuracy under different noise levels. For both the Minnesota and Bunny graphs, we fix the number of spatial samples per time step to \( m_t = 11 \) for regime 1 and \( m_t = 7 \) for regime 2. We then vary \( s \) from 2 to 15 in increments of 1 under different noise levels, with \( \sigma \in \{10^{-6}, 10^{-5}, 10^{-4}\} \). For each combination of the parameters, we perform 100 independent trials. The results are shown \Cref{fig_s_minnesota,fig_s_bunny}, respectively. In particular, when \( m_t = 7 \) and \( s = 2 \) in regime 2, reconstruction is unstable and even fails for the Bunny graph. However, increasing the number of time steps significantly improves reconstruction accuracy. Meanwhile, across all noise levels, the reconstruction error tends to decrease as \( s \) increases, with this trend being particularly pronounced in regime 2. This is consistent with the theoretical prediction in \cite[Proposition~4.1]{AHKLLV18}.

\subsubsection{Reconstruction without $U_k$} 
\begin{figure}[ht]
    \centering
    \begin{minipage}{0.24\linewidth}
        \centering
    \includegraphics[width=\linewidth]{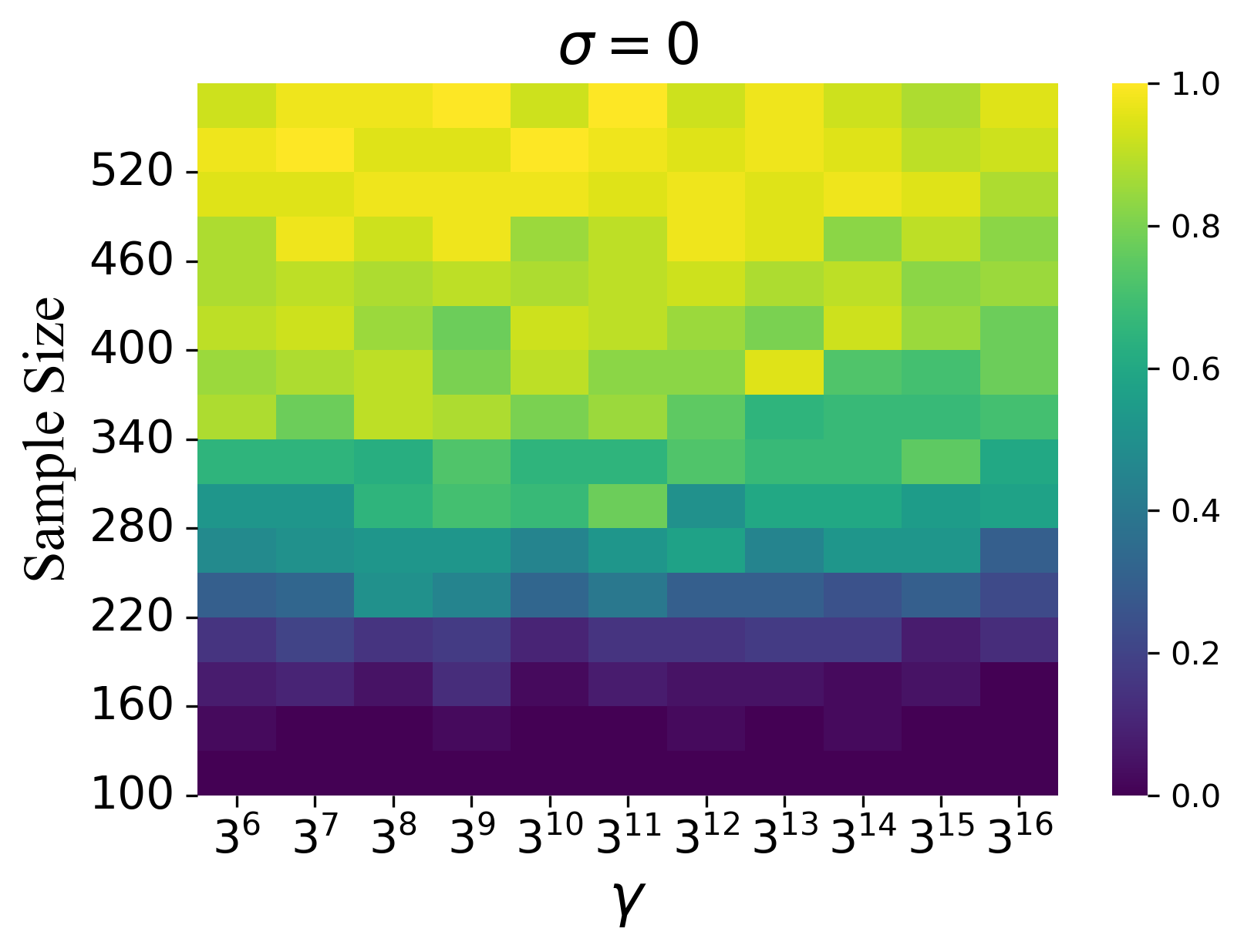}
    \end{minipage}
    \begin{minipage}{0.24\linewidth}
        \centering
    \includegraphics[width=\linewidth]{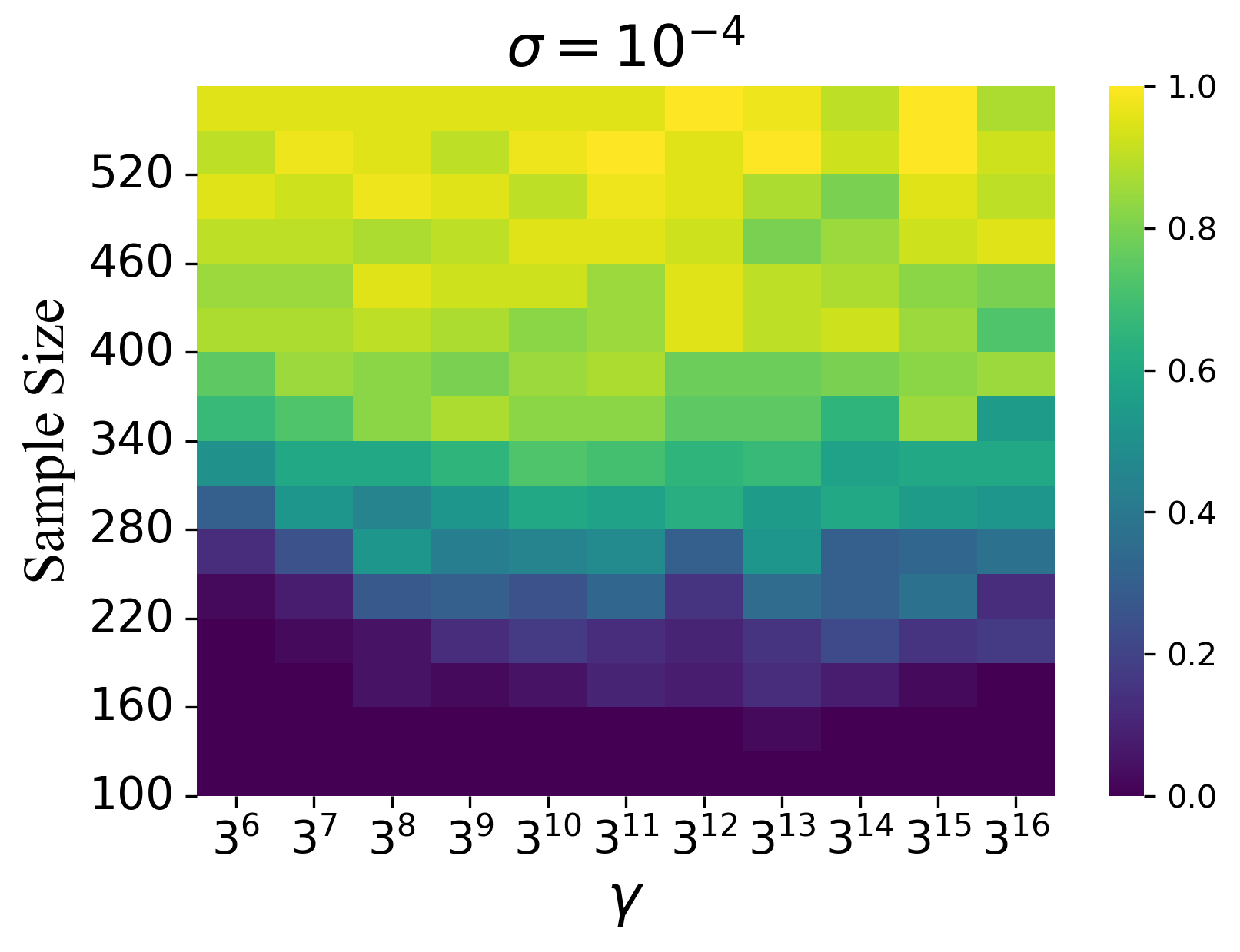}
    \end{minipage}
    \begin{minipage}{0.24\linewidth}
        \centering
    \includegraphics[width=\linewidth]{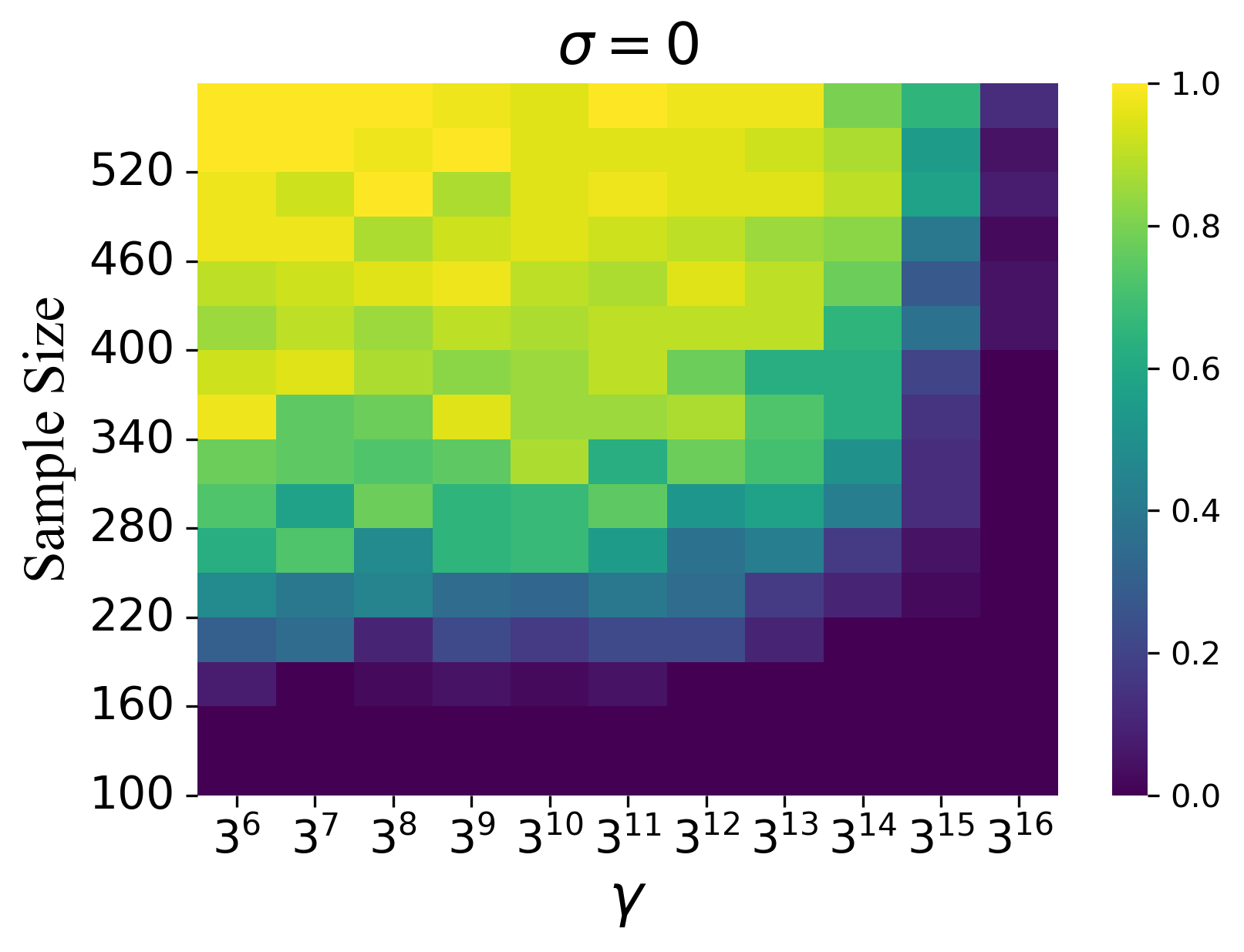}
    \end{minipage}
    \begin{minipage}{0.24\linewidth}
        \centering
    \includegraphics[width=\linewidth]{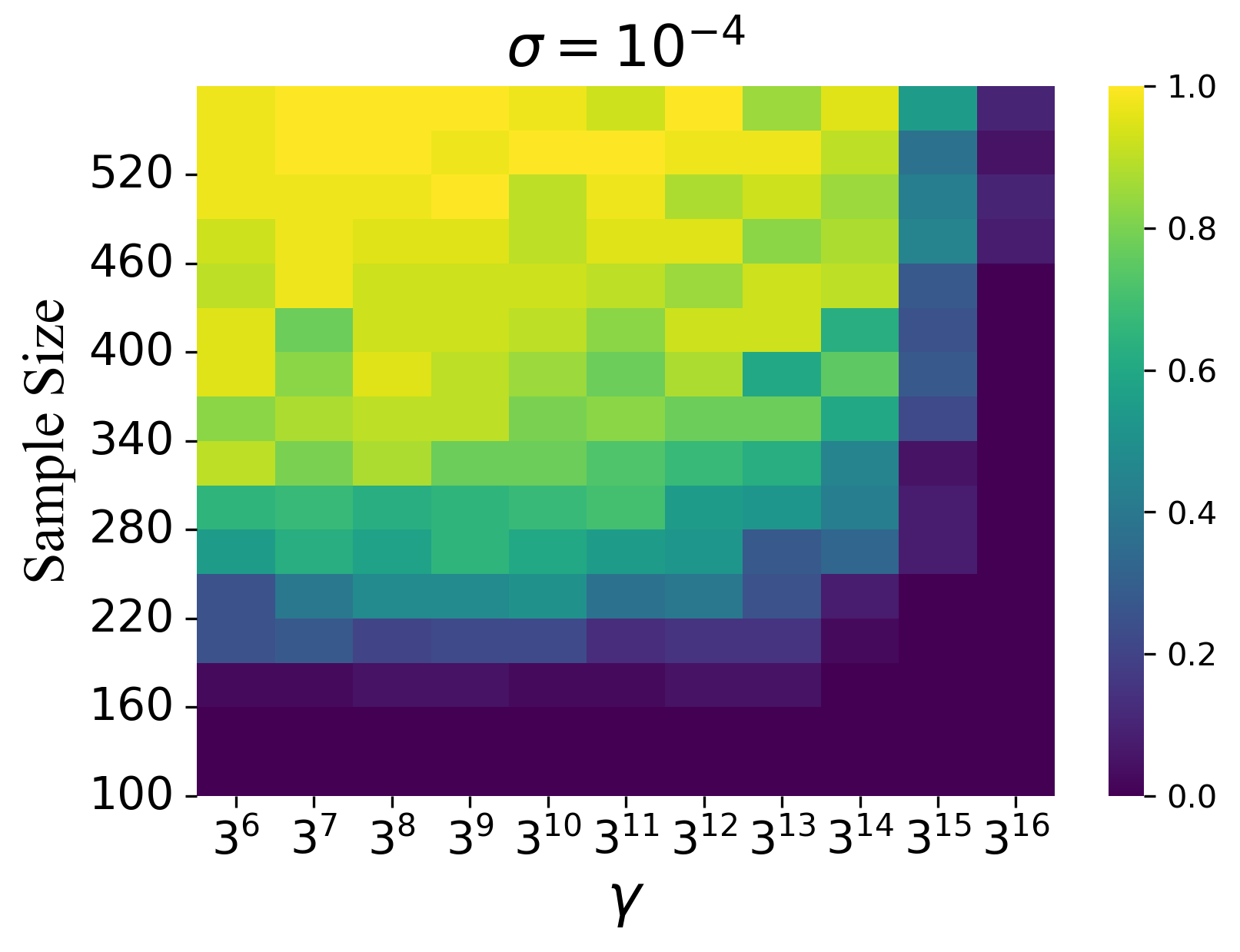}
    \end{minipage}
    
    \begin{minipage}{0.24\linewidth}
        \centering
    \includegraphics[width=\linewidth]{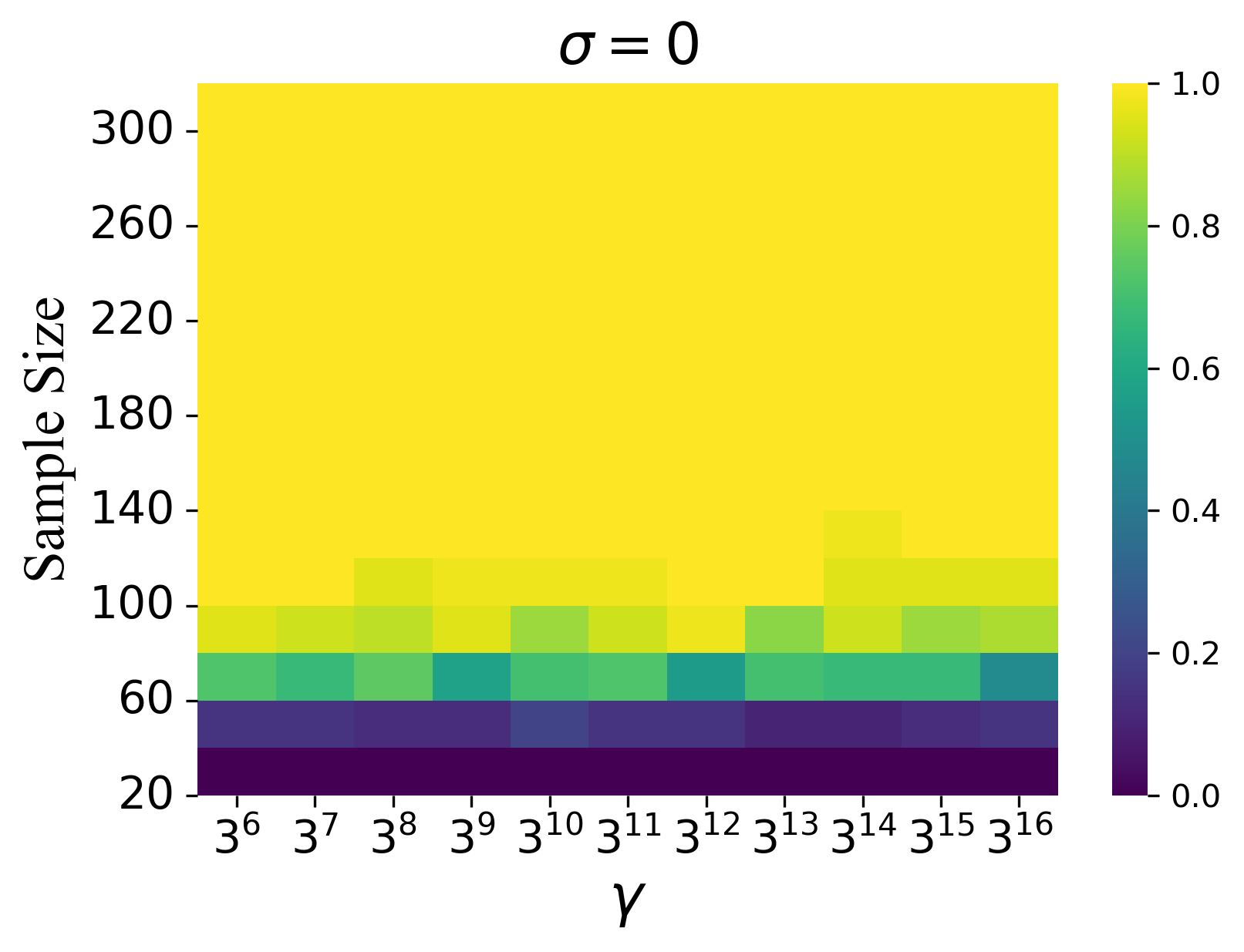}
    \end{minipage}
    \begin{minipage}{0.24\linewidth}
        \centering
    \includegraphics[width=\linewidth]{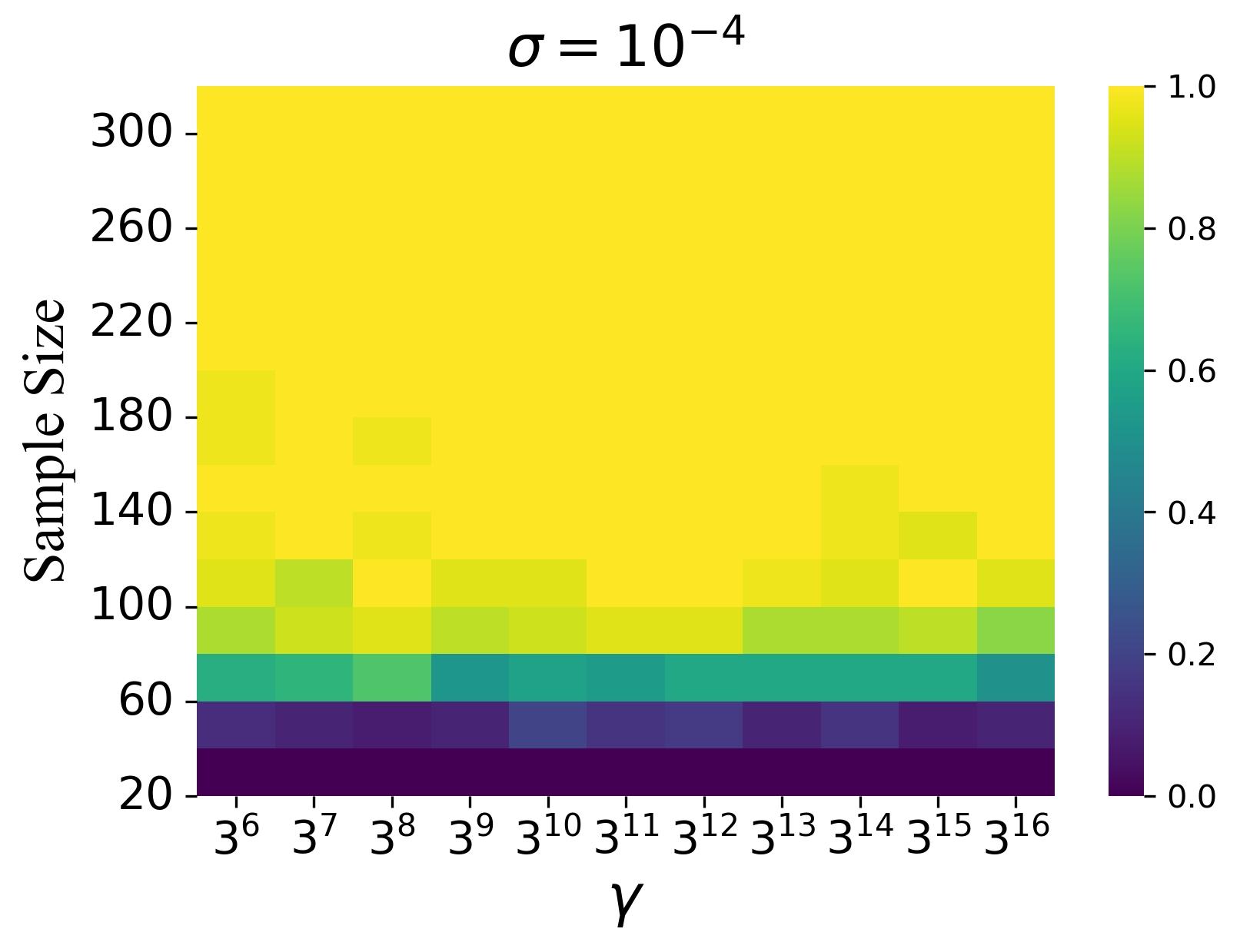}
    \end{minipage}
     \begin{minipage}{0.24\linewidth}
        \centering
    \includegraphics[width=\linewidth]{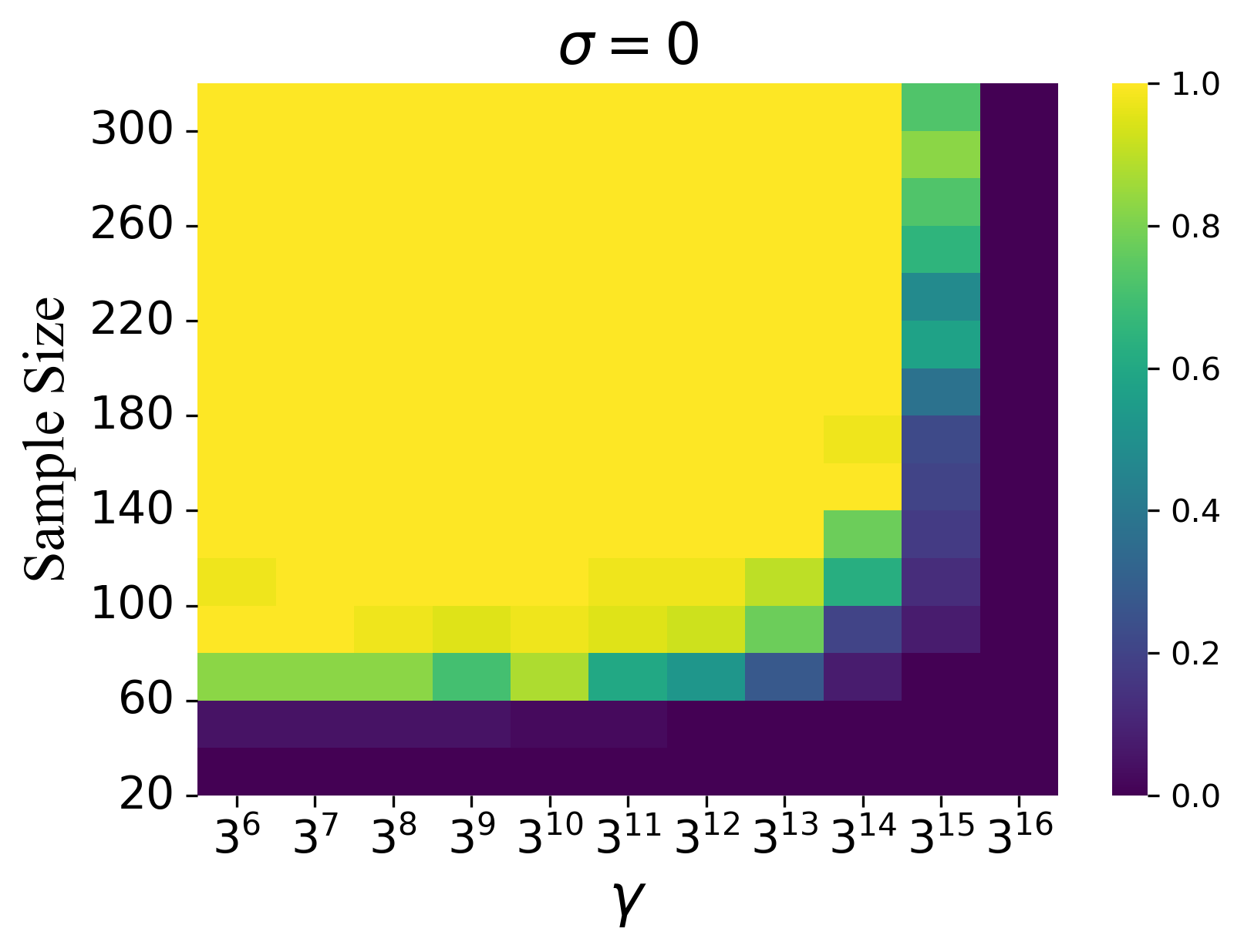}
    \end{minipage}
    \begin{minipage}{0.24\linewidth}
        \centering
    \includegraphics[width=\linewidth]{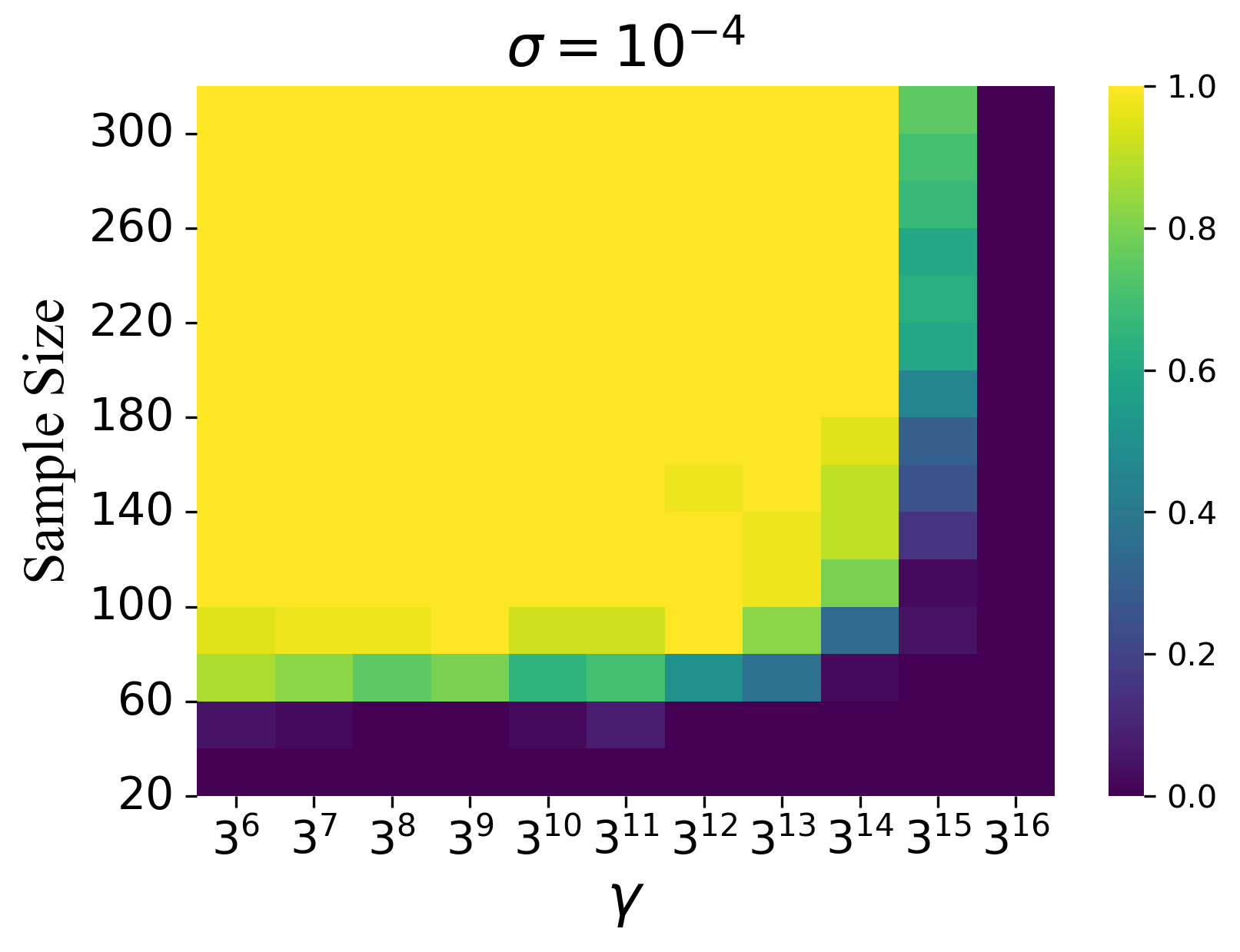}
    \end{minipage}
    \caption{Probability of \( \text{RE} < 0.05 \) for the Minnesota graph (first two columns) and the Bunny graph (last two columns). The first row corresponds to regime 1, and the second row corresponds to regime 2. The first and third column show the results without noise, while the second and fourth columns display the results with a noise level of $\sigma=10^{-4}$.}
  \label{fig3}
\end{figure}


When \( U_k \) is unknown, the vector \( \bw_{\text{aug}} \) is recovered using \eqref{recons2}. In this simulation, we set \( g(L) = L^4 \) and vary the sample size \( M \), the noise level \( \sigma \), and the parameter \( \gamma \). The values of \( \gamma \) are chosen from the set \( \{ 3^i \mid i \in \{6, 7, \dots, 16\} \} \), while \( \sigma \) varies over the set \(\{0, 10^{-4}\} \). For each combination of fixed parameters, we repeat the experiment 40 times. Since the differences in the results are not as significant as in the experiments with known \( U_k \), we plot heatmaps of the probability that the relative error is less than 0.05 to better illustrate the comparison, instead of using the boxplot as in the previous experiments. The results for the Minnesota and Bunny graphs are shown in  \Cref{fig3}, respectively.

Like the results in  \Cref{nonreg}, for both graphs, the experiment achieves better results with a much smaller sample size in regime 2 for certain scales of \( \gamma \) in both the noiseless and noisy cases. Additionally, we observe that the structure of the graphs has a significant influence on the reconstruction in these experiments. For instance, in regime 1, for both graphs, and for a fixed sample size \( M \), the results remain relatively stable within a certain range of \( \gamma \). However, the tolerance with respect to \( \gamma \) differs between these two graphs. The Minnesota graph exhibits higher tolerance with respect to \( \gamma \), as its performance remains stable even when \( \gamma \) exceeds \( 3^{14} \), while the Bunny graph starts to degrade. On the other hand, the Bunny graph shows higher tolerance with respect to noise. Specifically, when noise is added, the optimal \( \gamma \) shifts in both graphs, with the Minnesota graph requiring a larger \( \gamma \) to maintain similar reconstruction quality. Meanwhile, the Bunny graph remains more stable and less affected by noise.

\subsection{Experiments with Real-world Data}
In the experiments with real-world data, we did not evaluate the errors with respect to the initial value and the source term, as this was not feasible. Instead, we used the data from the first ten time instances to estimate the parameter \(\alpha\) and employed the approximated operator \(A=e^{-\alpha L}\) to complete the sampling-reconstruction process using the remaining data. From this, we obtained the approximated values for \(\bx_0\) and \(\bw\). We then used the relation \(\bx_{t+1} = A \bx_t + \bw\) to reconstruct the data and compared the recovered signal with the real signal at the unsampled points. 

In the sea surface temperature experiment, we employed a combination of global optimization (using PSO) and local optimization (via fmincon with SQP) to find the optimal parameter \(\alpha\). However, in the COVID-19 experiment, this approach was ineffective due to the ill-conditioned nature of the dataset, which prevented us from finding an optimal parameter with a low error. As a result, we used a greedy algorithm to find the optimal value of \(\alpha\). Additionally, in the sea surface temperature experiment, we set the bandwidth to 10, while in the COVID-19 experiment, we set it to 38. This choice ensures that at least 90 percent of the signal energy is concentrated within the selected bandwidth.
\begin{figure}[ht]
    \centering
    \begin{minipage}{0.32\linewidth}
        \centering
    \includegraphics[width=\linewidth]{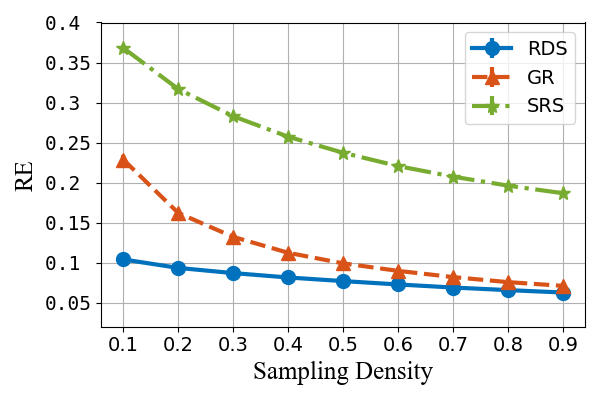}
    \end{minipage}
    \begin{minipage}{0.32\linewidth}
        \centering
\includegraphics[width=\linewidth]{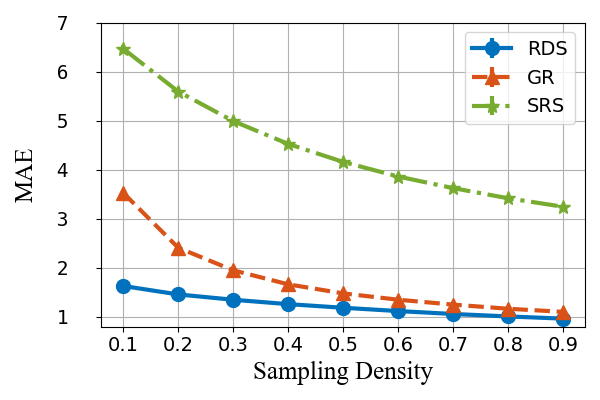}
    \end{minipage}
    \begin{minipage}{0.32\linewidth}
        \centering
    \includegraphics[width=\linewidth]{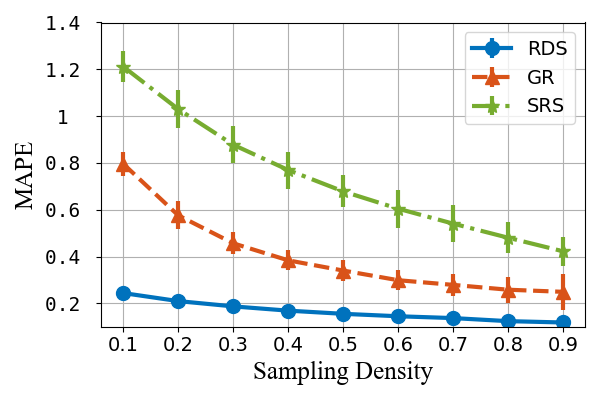}
    \end{minipage}

     \begin{minipage}{0.323\linewidth}
        \centering        \includegraphics[width=\linewidth]{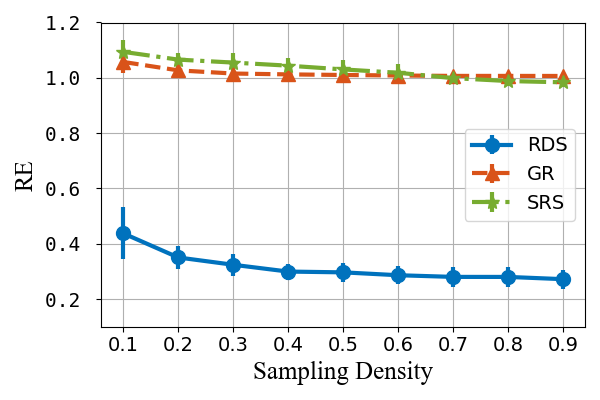}
    \end{minipage}
    \begin{minipage}{0.323\linewidth}
        \centering
      \includegraphics[width=\linewidth]{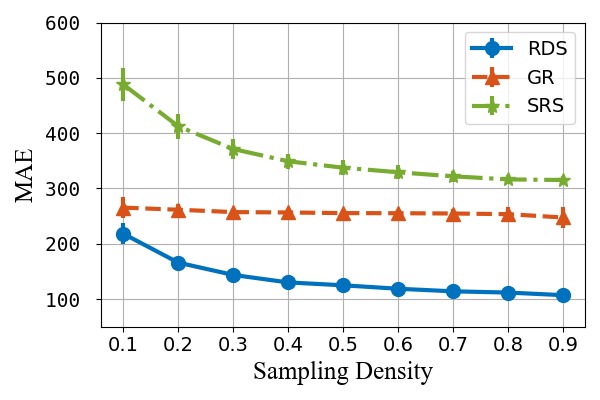}
    \end{minipage}
    \begin{minipage}{0.323\linewidth}
        \centering
   \includegraphics[width=\linewidth]{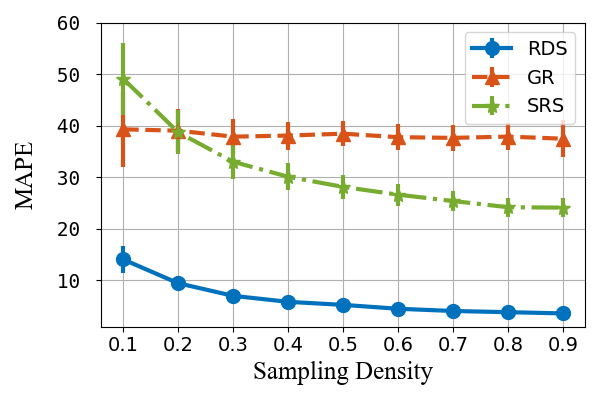}
    \end{minipage}
    \caption{Comparison of RDS, GR, and SRS in reconstruction across different sampling densities on Sea Surface Temperature dataset (first row) and Global COVID-19 dataset (second row).}
    \label{fig5}
\end{figure}

There are other methods applied to recover the same dataset \cite{9730033}, but they rely on different assumptions. Therefore, we limit our comparison to the methods based on Graph Regularization (GR) \cite{6736922} and Static Random Sampling (SRS) \cite{puy2018random}, which, like ours, assume the signal is bandlimited. Due to the complex nature of real data, experiments in regime 1 could result in significant errors. Therefore, we conducted experiments only in regime 2. We set the sampling rate per time instance to range from 0.1 to 0.9 in steps of 0.1. In other words, the total space–time samples are $[0.1:0.1:0.9]\times ns$ where $n$ is the number of graph nodes and $s$ is the total sampling time. For each dataset, we conducted 100 experiments for each sampling density, computed three different evaluation metrics, and plotted the results with error bars, as shown in  \Cref{fig5}. The results demonstrate that our method (RDS) achieves superior performance compared to other approaches across various evaluation metrics.


\section*{Acknowledgment}
This work is partially supported by NSF grant DMS-2208030.

\appendix
\section{Proofs}
\begin{proof}[The proof of \Cref{embedding}]
\begin{eqnarray*}
\|\pi_{A,s}(\bv)\|^2
&=&
\left(\begin{bmatrix}
        U_k^\top&\bf{0}\\
        \bf{0}&U_k^\top
    \end{bmatrix}\bv\right)^{\top}\left(\sum_{t=0}^{s-1}
\begin{bmatrix}
    \Lambda_k^{t}  \\
    \bar{\Lambda}_k^{t}  
\end{bmatrix}
    U_k^{\top} 
  U_k
\begin{bmatrix}
    \Lambda_k^{t} & \bar{\Lambda}_k^{t}
\end{bmatrix}\right)\left(\begin{bmatrix}
        U_k^\top&\bf{0}\\
        \bf{0}&U_k^\top
    \end{bmatrix}\bv\right)\\
&=&
\left(\begin{bmatrix}
        U_k^\top&\bf{0}\\
        \bf{0}&U_k^\top
    \end{bmatrix}\bv\right)^{\top}\left(
\sum_{t=0}^{s-1}
\underbrace{
\begin{bmatrix}
    \Lambda_k^{2t} & \Lambda_k^{t}\bar{\Lambda}_k^{t}\\
    \Lambda_k^{t}\bar{\Lambda}_k^{t} & \bar{\Lambda}_k^{t}\bar{\Lambda}_k^{t}
\end{bmatrix}}_{:=Y_t}
\right)\left(\begin{bmatrix}
        U_k^\top&\bf{0}\\
        \bf{0}&U_k^\top
    \end{bmatrix}\bv\right).\\
\end{eqnarray*}
Observe that for each $t,$ $Y_t$ is a self-adjoint, positive semidefinite matrix with $\Lambda^0=I$, $\bar{\Lambda}^0=\bf{0}$ and for $t\ge 1$
\[
 \Lambda_k^{2t}=\begin{bmatrix}
        \lambda_1^{2t}&&&\\
         &\lambda_2^{2t}&&\\
      & &\ddots&\\
       & &&\lambda_k^{2t} 
    \end{bmatrix}, 
\Lambda_k^{t}\bar{\Lambda}_k^{t}=
\begin{bmatrix}
        \lambda_1^{t}(\sum\limits_{l=0}^{t-1}\lambda_1^l) & & &\\
        &  \lambda_2^{t}(\sum\limits_{l=0}^{t-1}\lambda_2^l) & & \\
      &&\ddots& \\
         && &\lambda_k^{t}(\sum\limits_{l=0}^{t-1}\lambda_k^l)
    \end{bmatrix}
\]
and 
\[
\bar{\Lambda}_k^{t}\bar{\Lambda}_k^{t}=
\begin{bmatrix}
        (\sum\limits_{l=0}^{t-1}\lambda_1^l)^2 & && \\
        & (\sum\limits_{l=0}^{t-1}\lambda_2^l)^2 & & \\
      & &\ddots& \\
        & & &(\sum\limits_{l=0}^{t-1}\lambda_k^l)^2 
    \end{bmatrix}. 
\]
Therefore, $Y:=\sum_{t=0}^{s-1} Y_t$ forms a self-adjoint, positive semidefinite matrix, for which we seek to determine the largest and smallest eigenvalues. To solve this, we notice that there exists some permutation matrix $P$ such that $\widetilde{Y}:=P^\top YP$ is a block diagonal matrix. More specifically, $\widetilde{Y}$ equals 
\[
\begin{bmatrix}
\widetilde{Y}_{1,1}& && \\
&\widetilde{Y}_{2,2}& & \\
& &\ddots& \\
 && &\widetilde{Y}_{k,k}
\end{bmatrix}
\]
where 
\[
\widetilde{Y}_{j,j}=
\begin{bmatrix}
\sum\limits_{t=0}^{s-1}\lambda_j^{2t} & \sum\limits_{t=0}^{s-1}(\sum\limits_{l=0}^{t-1}\lambda_j^{l+t})\\
\sum\limits_{t=0}^{s-1}(\sum\limits_{l=0}^{t-1}\lambda_j^{l+t}) & \sum\limits_{t=0}^{s-1}(\sum\limits_{l=0}^{t-1}\lambda_j^l)^2 
\end{bmatrix} 
\]
for $j\in [k].$

This  implies that 
\begin{equation}\label{boundc}
c_{A,k,s}:=\lambda_{min}\left(Y\right)=\lambda_{min}\left(\widetilde{Y}\right)=\min\limits_{1\le j\le k} \lambda\left(\widetilde{Y}_{j,j}\right)
\end{equation}
and
\begin{equation}\label{boundC}
C_{A,k,s}:=\lambda_{max}\left(Y\right)=\lambda_{max}\left(\widetilde{Y}\right)=\max\limits_{1\le j\le k} \lambda\left(\widetilde{Y}_{j,j}\right).
\end{equation}
Next, let's show $c_{A,k,s}>0$. By Cauchy–Schwarz inequality, we have 
\[
\left(\sum\limits_{t=0}^{s-1}(\sum\limits_{l=0}^{t-1}\lambda_j^{l+t})\right)^2=
\left(\sum\limits_{t=0}^{s-1}(\lambda_j^t\sum\limits_{l=0}^{t-1}\lambda_j^{l})\right)^2\le \left(\sum\limits_{t=0}^{s-1}\lambda_j^{2t}\right)\left(\sum\limits_{t=0}^{s-1}(\sum\limits_{l=0}^{t-1}\lambda_j^l)^2 \right),
\]
where the equality holds if and only if there exists some constant $c$ such that $\lambda_j^t=c\sum\limits_{l=0}^{t-1}\lambda_j^l$ for all $t,$ which is impossible. Therefore, $\text{det}(\widetilde{Y}_{j,j})>0$ for all $j\in [k],$ implying that $C_{A,k,s}\ge c_{A,k,s}>0.$

As a result, we can derive the following inequality
\[
c_{A,k,s}\|\bv\|^2=c_{A,k,s}\left\|\begin{bmatrix}
        U_k^\top&\bf{0}\\
        \bf{0}&U_k^\top
    \end{bmatrix}\bv\right\|^2 \le \|\pi_{A,s}(\bv)\|^2\le C_{A,k,s}\left\|\begin{bmatrix}
        U_k^\top&\bf{0}\\
        \bf{0}&U_k^\top
    \end{bmatrix}\bv\right\|^2=C_{A,k,s}\|\bv\|^2.
\]
\end{proof}

\begin{proof}[The proof of \Cref{nuestimate}]
    Notice that 
\begin{eqnarray*}
    \frac{1}{\textbf{p}^{(2)}_t(\ell)}\|\delta_{\ell}^\top U_k
    \begin{bmatrix}
        \Lambda_k^t & \bar{\Lambda}_k^t
    \end{bmatrix}\|_2^2 &\le& \frac{\|\delta_{\ell}^\top U_k\|^2}{\textbf{p}^{(2)}_t(\ell)}
    \|\begin{bmatrix}
        \Lambda_k^t & \bar{\Lambda}_k^t
    \end{bmatrix}\|_2^2\\
    &=&\frac{\|\delta_{\ell}^\top U_k\|^2}{\textbf{p}^{(2)}_t(\ell)}
    \left\|\begin{bmatrix}
        \Lambda_k^{2t} &\Lambda_k^t\bar{\Lambda}_k^t\\
        \Lambda_k^t\bar{\Lambda}_k^t & \bar{\Lambda}_k^t\bar{\Lambda}_k^t
    \end{bmatrix}\right\|_2\\
    &=&  \frac{\|U_k^\top\delta_\ell\|^2}{\mathbf{p}^{(2)}_t(\ell)} \max_{1\leq j\leq k}\left\|\begin{bmatrix}
 \lambda_j^{2t}   &   \sum\limits_{l=0}^{t-1}\lambda_j^{l+t}  \\
 \sum\limits_{l=0}^{t-1}\lambda_j^{l+t}  &   (\sum\limits_{l=0}^{t-1}\lambda_j^l)^2 \end{bmatrix}\right\|_2\\
 &\le&\max_{1\leq i\leq n}\frac{\|U_k^\top\delta_i\|^2}{\mathbf{p}^{(2)}_t(i)} \max_{1\leq j\leq k} (\lambda_j^{2t}+(\sum_{i=0}^{t-1}\lambda_j^t)^2).
\end{eqnarray*}

Additionally, since the matrix $\begin{bmatrix}
    \Lambda_k^{t} \\
    \bar{\Lambda}_k^{t}
\end{bmatrix}
U_k^{\top}\delta_{\ell}\delta_{\ell}^{\top}U_k
\begin{bmatrix}
    \Lambda_k^{t} & \bar{\Lambda}_k^{t}\\
\end{bmatrix}$ is symmetric and positive semidefinite for all $t\in\{0,1,\cdots,s-1\},$ we have $$\left\|
   \delta_\ell^{\top}U_k
\begin{bmatrix}
    \Lambda_k^{t} & \bar{\Lambda}_k^{t}\\
\end{bmatrix}\right\|_2^2
=\left\|\begin{bmatrix}
    \Lambda_k^{t} \\
    \bar{\Lambda}_k^{t}
\end{bmatrix}
U_k^{\top}\delta_{\ell}\delta_{\ell}^{\top}U_k
\begin{bmatrix}
    \Lambda_k^{t} & \bar{\Lambda}_k^{t}\\
\end{bmatrix}\right\|_2
\le\left\|\sum\limits_{t=0}^{s-1}\begin{bmatrix}
    \Lambda_k^{t} \\
    \bar{\Lambda}_k^{t}
\end{bmatrix}
U_k^{\top}\delta_{\ell}\delta_{\ell}^{\top}U_k
\begin{bmatrix}
    \Lambda_k^{t} & \bar{\Lambda}_k^{t}\\
\end{bmatrix}\right\|_2.$$ 
Thus, when $\mathbf{P}^{(j)}\ (j=1,2)$ are uniform distributions, $\nu^{(1)}\ge\nu^{(2)}(t)$ for all $t=0,1,\cdots,s-1.$
\end{proof}

\begin{proof}[The proof of \Cref{thm:smp_cmp}]
We first present the proof for regime 2. Recall that $\widetilde{U}_k=\begin{bmatrix}
    U_k &\textbf{0}\\
    \textbf{0} & U_k
\end{bmatrix}.$ For $\bw_{\text{aug}}\in\text{span}(\widetilde{U}_k),$ then 
\[
\begin{bmatrix}
   \bx_0\\
    \bx_1\\
    \vdots\\
    \bx_{s-2}\\
    \bx_{s-1}
\end{bmatrix}=\pi_{A,s}(\bw_{\text{aug}})=\begin{bmatrix}
U_k\begin{bmatrix}
    I &\bf{0}
\end{bmatrix}\\
U_k
\begin{bmatrix}
    \Lambda_k & I
\end{bmatrix}\\
\vdots\\
    U_k 
\begin{bmatrix}
    \Lambda_k^{s-1} & \bar{\Lambda}_k^{s-1}
\end{bmatrix}
\end{bmatrix}
\hat{\bw}_{\text{aug}}
\]
where $\hat{\bw}_{\text{aug}}=\begin{bmatrix}
    \hat{\bx}_0\\
    \hat{\bw}
\end{bmatrix}.$ Next, we obtain
\begin{eqnarray*}
\mathcal{W}\mathcal{P}\mathcal{S}
   \begin{bmatrix}
   \bx_0\\
    \bx_1\\
    \vdots\\
    \bx_{s-2}\\
    \bx_{s-1}
\end{bmatrix}
=\underbrace{\begin{bmatrix}
    M_0 P_{\Omega_0}^{-\frac{1}{2}}S_0U_k
    \begin{bmatrix}
        I &\bf{0}
    \end{bmatrix}\\
M_1 P_{\Omega_1}^{-\frac{1}{2}}S_1 U_k\begin{bmatrix}
       \Lambda_{k}^{1}&\bar{\Lambda}_k^{1}
   \end{bmatrix}\\
   \vdots\\
   M_{s-2} P_{\Omega_{s-2}}^{-\frac{1}{2}}S_{s-2} U_k\begin{bmatrix}
       \Lambda_{k}^{s-2}&\bar{\Lambda}_k^{s-2}
   \end{bmatrix}\\
   M_{s-1} P_{\Omega_{s-1}}^{-\frac{1}{2}}S_{s-1} U_k\begin{bmatrix}
       \Lambda_{k}^{s-1}&\bar{\Lambda}_k^{s-1}
   \end{bmatrix}
\end{bmatrix}}_B  
   \begin{bmatrix}
       \hat{\bx}_0\\
       \hat{\bw}
   \end{bmatrix}.
\end{eqnarray*}

$B$ is defined as above, then we get
\begin{eqnarray*}
    B^{\top}B
    &=& \sum_{t=0}^{s-1}\begin{bmatrix}
        \Lambda_k^t\\
        \bar{\Lambda}_k^{t}
        \end{bmatrix}U_k^{\top}S_t^{\top}P_{\Omega_t}^{\top}M_{t}^{\top}M_{t}P_{\Omega_t}S_t U_k\begin{bmatrix}
       \Lambda_{k}^{t-1}&\bar{\Lambda}_k^{t-1}
   \end{bmatrix}\\
    &=& 
    \sum_{t=0}^{s-1}\begin{bmatrix}
        \Lambda_k^t\\
        \bar{\Lambda}_k^{t}
        \end{bmatrix}U_k^{\top} \sum_{j=1}^{m_t}\frac{\delta_{\omega_{j,t}}\delta_{\omega_{j,t}}^{\top}}{m_t\mathbf{p}_t(\omega_{j,t})}
        U_k\begin{bmatrix}
       \Lambda_{k}^{t-1}&\bar{\Lambda}_k^{t-1}\end{bmatrix}\\
&=&
 \sum_{t=0}^{s-1}\sum_{j=1}^{m_t}
\begin{bmatrix}
    \Lambda_k^{t} \\
    \bar{\Lambda}_k^{t}
\end{bmatrix}
\frac{U_k^{\top}\delta_{\omega_{j,t}}\delta_{\omega_{j,t}}^{\top}U_k}{m_t\mathbf{p}_t(\omega_{j,t})}
\begin{bmatrix}
    \Lambda_k^{t} & \bar{\Lambda}_k^{t}\\
\end{bmatrix}.
\end{eqnarray*}

Let us define
\begin{eqnarray*}
X=B^{\top}B:=\sum_{t=0}^{s-1}\sum_{j=1}^{m_t} X_{t,\omega_{j,t}}
\end{eqnarray*}
where 
\[X_{t,\omega_{j,t}}=\frac{1}{{m_{t}{\mathbf{p}}_{t}( \omega_{j,t}}) }\begin{bmatrix}
    \Lambda_k^{t} \\
    \bar{\Lambda}_k^{t}
\end{bmatrix}
 U_k^{\top}\delta_{\omega_{j,t}}\delta_{\omega_{j,t}}^{\top}U_k
\begin{bmatrix}
    \Lambda_k^{t} & \bar{\Lambda}_k^{t}\\
\end{bmatrix}.\]

Since all $X_{t,\omega_{j,t}}$ are self-adjoint, positive semidefinite matrices, we have 
\begin{eqnarray*}
    \lambda_{max}(X_{t,\omega_{j,t}})
    &=&\|X_{t,\omega_{j,t}}\|_2\le\max_{0\le t\le s-1,\ \omega_{j,t}\in\Omega_t} \|X_{t, \omega_{j,t}}\|_2\\
    &=& \frac{1}{m_t}\max_{0\le t\le s-1,\ \omega_{j,t}\in\Omega_t}
    \left\{\frac{1}{\textbf{p}_t(\omega_{j,t})}\left\|
   \delta_{\omega_{j,t}}^{\top}U_k
\begin{bmatrix}
    \Lambda_k^{t} & \bar{\Lambda}_k^{t}\\
\end{bmatrix}\right\|_2^2\right\}\leq\max_{0\le t\le s-1}\frac{\nu^{(2)}(t)}{m_t}:=\nu.
\end{eqnarray*}

Since each node $\omega_{j,t}$ in the sampling set $\Omega_t$ is randomly and independently selected from $[n]$ according to probability distribution $\mathbf{p}_t$, we have
\begin{eqnarray*}
    \mathbb{E}(X_{t,\omega_{j,t}}) &=& \frac{1}{m_t}
\begin{bmatrix}
    \Lambda_k^{t}  \\
    \bar{\Lambda}_k^{t} 
\end{bmatrix}
 \mathbb{E}\left( 
\frac{U_k^{\top}\delta_{\omega_{j,t}}\delta_{\omega_{j,t}}^{\top}U_k}{\textbf{p}_t(\omega_{j,t})}  \right)
\begin{bmatrix}
    \Lambda_k^{t} & \bar{\Lambda}_k^{t}
\end{bmatrix}\\
&=& 
\frac{1}{m_t}
\begin{bmatrix}
    \Lambda_k^{t}  \\
    \bar{\Lambda}_k^{t} 
\end{bmatrix}\left(\sum_{\ell=1}^N\textbf{p}_t(\ell) 
\frac{U_k^{\top}\delta_{\omega_{j,t}}\delta_{\omega_{j,t}}^{\top} U_k}{\textbf{p}_t(\ell)}  \right)
\begin{bmatrix}
    \Lambda_k^{t} & \bar{\Lambda}_k^{t}
\end{bmatrix}\\
&=&
\frac{1}{m_t}
\begin{bmatrix}
    \Lambda_k^{t}  \\
    \bar{\Lambda}_k^{t} 
\end{bmatrix}\begin{bmatrix}
    \Lambda_k^{t} & \bar{\Lambda}_k^{t}
\end{bmatrix}\\
&=&
\frac{1}{m_t}\underbrace{
\begin{bmatrix}
    \Lambda_k^{2t}  & \Lambda_k^{t}\bar{\Lambda}_k^{t}\\
  \Lambda_k^{t}\bar{\Lambda}_k^{t}& (\bar{\Lambda}_k^{t})^2
\end{bmatrix}}_{Y_t}\\
\end{eqnarray*}
for all $X_{t, \omega_{j,t}}$ for some fixed $t$.

According to \eqref{boundc} and \eqref{boundC}, we obtain
\[
\lambda_{min}(\sum_{t,j}\mathbb{E}(X_{t, \omega_{j,t}}))=\lambda_{min}\left(\sum\limits_{t=0}^{s-1} Y_t\right)=c_{A,k,s}
\]
and
\[
\lambda_{max}(\sum_{t,j}\mathbb{E}(X_{t, \omega_{j,t}}))=\lambda_{max}\left(\sum\limits_{t=0}^{s-1} Y_t\right)=C_{A,k,s}.
\]
According to the concentrated inequality in \cite[Theorem 1.1]{tropp2012user}, we thus have
\begin{eqnarray*}
    \mathbb{P}\{\lambda_{min}(X)\le (1-\delta)c_{A,k,s}\}
    \le 
    2k\left[\frac{e^{-\delta}}{(1-\delta)^{1-\delta}}\right]^{\frac{ c_{A,k,s}}{\nu}}
    \le 2k\ \exp (-\frac{\delta^2c_{A,k,s}}{3\nu}),
\end{eqnarray*}

\begin{eqnarray*}
    \mathbb{P}\{\lambda_{\max}(X)\ge (1+\delta)C_{A,k,s}\}
    \le
    2k\left[\frac{e^{\delta}}{(1+\delta)^{1+\delta}}\right]^{\frac{C_{A,k,s}}{\nu}}
    \le 2k\ \exp (-\frac{\delta^2 C_{A,k,s}}{3\nu})
    \le 2k\ \exp (-\frac{\delta^2 c_{A,k,s}}{3\nu}).
\end{eqnarray*}
Obviously, $2k\ \exp (-\frac{\delta^2c_{A,k,s}}{3\nu})\le\frac{\epsilon}{2}$ always holds if

\begin{equation}\label{nu}
    \frac{1}{\nu}\ge \frac{3}{c_{A,k,s}\delta^2}\log\frac{4k}{\epsilon}.
\end{equation}

Notice that 
\begin{equation*}
    m_t\ge \frac{3\nu^{(2)}(t)}{c_{A,k,s}\delta^2}\log\frac{4k}{\epsilon}
\end{equation*}
for all  $t=0,\cdots,s-1$ can ensure \eqref{nu} holds.

Finally, we get that
\[
(1-\delta)c_{A,k,s}\left\|
\bw_{\text{aug}}\right\|_2^2
\le
\left\|\mathcal{W}\mathcal{P}\mathcal{S}
    \begin{bmatrix}
    \bx_0\\
    \vdots\\
    \bx_{s-1}
\end{bmatrix}\right\|^2_2
=\left\|\mathcal{W}\mathcal{P}\mathcal{S}\pi_{A,s}(\bw_{\text{aug}})\right\|_2^2\le(1+\delta)C_{A,k,s}\left\|\bw_{\text{aug}}\right\|_2^2.
\]

For regime 1, $\Omega_{t}:=\Omega$, $m_t:=m$ for $t=0,\cdots,s-1$, and $\Omega$ is determined by $\mathbf{p}_0$.  Thus
\begin{eqnarray*}
X=\sum_{t=0}^{s-1}\sum_{j=1}^m X_{t,\omega_{j}}=\sum_{j=1}^m X_j 
\end{eqnarray*}
where
\[X_j=\sum_{t=0}^{s-1}X_{t,\omega_j}=\frac{1}{{m{\mathbf{p}_0}(\omega_j)} }
\sum_{t=0}^{s-1}\begin{bmatrix}
    \Lambda_k^{t} \\
    \bar{\Lambda}_k^{t}
\end{bmatrix}
U_k^{\top}\delta_{\omega_{j}}\delta_{\omega_{j}}^{\top}U_k
\begin{bmatrix}
    \Lambda_k^{t} & \bar{\Lambda}_k^{t}\\
\end{bmatrix}.\]

Since all $X_j$ are self-adjoint, positive semidefinite matrices, we have 
\begin{eqnarray*}
    \lambda_{max}(X_j)
    &=&\|X_{j}\|_2\le\max_{\omega_{j}\in\Omega} \|X_{j}\|_2\\
    &=& \frac{1}{m}\max_{\omega_{j}\in\Omega}
    \left\{\frac{1}{\mathbf{p}_0(\omega_{j})}\left\|\sum_{t=0}^{s-1}\begin{bmatrix}
    \Lambda_k^{t} \\
    \bar{\Lambda}_k^{t}
\end{bmatrix}
U_k^{\top}\delta_{\omega_{j}}\delta_{\omega_{j}}^{\top}U_k
\begin{bmatrix}
    \Lambda_k^{t} & \bar{\Lambda}_k^{t}\\
\end{bmatrix}\right\|_2\right\}\leq\frac{\nu^{(1)}}{m}.
\end{eqnarray*}

Additionally, we have
\begin{eqnarray*}
    \mathbb{E}(X_j) &=& \sum_{t=0}^{s-1}\frac{1}{m}
\begin{bmatrix}
    \Lambda_k^{t}  \\
    \bar{\Lambda}_k^{t} 
\end{bmatrix}
 \mathbb{E}\left( 
\frac{U_k^{\top}\delta_{\omega_{j}}\delta_{\omega_{j}}^{\top}U_k}{{\mathbf{p}}(\omega_{j})}  \right)
\begin{bmatrix}
    \Lambda_k^{t} & \bar{\Lambda}_k^{t}
\end{bmatrix}\\
&=& 
\frac{1}{m}
\begin{bmatrix}
    \Lambda_k^{t}  \\
    \bar{\Lambda}_k^{t} 
\end{bmatrix}\left(\sum_{\ell=1}^N{\mathbf{p}}(\ell) 
\frac{U_k^{\top}\delta_{\omega_{j}}\delta_{\omega_{j}}^{\top} U_k}{\mathbf{p}(\ell)}  \right)
\begin{bmatrix}
    \Lambda_k^{t} & \bar{\Lambda}_k^{t}
\end{bmatrix}=
\frac{1}{m}\sum_{t=0}^{s-1}\underbrace{
\begin{bmatrix}
    \Lambda_k^{2t}  & \Lambda_k^{t}\bar{\Lambda}_k^{t}\\
  \Lambda_k^{t}\bar{\Lambda}_k^{t}& (\bar{\Lambda}_k^{t})^2
\end{bmatrix}}_{Y_t}\\
\end{eqnarray*}
for all $X_{j}.$ The rest of the proof follows similarly to the proof of regime 2.

\end{proof}

\begin{proof}[The proof of \Cref{alg1}]
\begin{enumerate}
\item
By optimality of $\bw_{\text{aug}}^*,$ we obtain 
\begin{equation*}
    \|\mathcal{W}\mathcal{P}(\mathcal{S}\pi_{A,s}(\bw_{\text{aug}}^*)-\bz)\|_2\le \|\mathcal{W}\mathcal{P}(\mathcal{S}\pi_{A,s}(\overline{\bw_{\text{aug}}})-\bz)\|_2
\end{equation*}
for any $\overline{\bw_{\text{aug}}}\in\text{span}(\widetilde{U}_k).$ Particularly, if $\overline{\bw_{\text{aug}}}= {\bw_{\text{aug}}}=\begin{bmatrix}
    \bx_0\\
    \bw
\end{bmatrix},$ then we have
\begin{equation}\label{errorbound1}
\|\mathcal{W}\mathcal{P}(\mathcal{S}\pi_{A,s}(\bw_{\text{aug}}^*)-\bz)\|_2\le \|\mathcal{W}\mathcal{P}(\mathcal{S}\pi_{A,s}(\bw_{\text{aug}})-\bz)\|_2=\|\mathcal{W}\mathcal{P}\be\|_2.
\end{equation}
On the other hand, we can obtain
\begin{eqnarray}\label{errorbound2}
    \|\mathcal{W}\mathcal{P}(\mathcal{S}\pi_{A,s}(\bw_{\text{aug}}^*)-\bz)\|_2
    &=&\|\mathcal{W}\mathcal{P}\mathcal{S}\pi_{A,s}(\bw_{\text{aug}}^*)-\mathcal{W}\mathcal{P}\mathcal{S}\pi_{A,s}(\bw_{\text{aug}})-\mathcal{W}\mathcal{P}\be\|_2\nonumber\\
    &\ge&
    \|\mathcal{W}\mathcal{P}\mathcal{S}\pi_{A,s}(\bw_{\text{aug}}^*-\bw_{\text{aug}})\|_2-\|\mathcal{W}\mathcal{P}\be\|_2\nonumber\\
    &\ge&
    \sqrt{(1-\delta)c_{A,k,s}}\left\|
\bw_{\text{aug}}^*-\bw_{\text{aug}}\right\|_2-\|\mathcal{W}\mathcal{P}\be\|_2.
\end{eqnarray}
Combing \eqref{errorbound1} and \eqref{errorbound2}, the following inequality holds
\begin{equation*}
    \left\|
\bw_{\text{aug}}^*-\bw_{\text{aug}}\right\|_2\le \frac{2}{\sqrt{(1-\delta)c_{A,k,s}}}\|\mathcal{W}\mathcal{P}\be\|_2
\end{equation*}
with probability at least $1-\epsilon$.
\item 
If we choose $\widetilde{\be}=\mathcal{S}\pi_{A,s}(\widetilde{\bw_{\text{aug}}})$ for some $\widetilde{\bw_{\text{aug}}}\in\text{span}(\widetilde{U}_k),$ then the space-time samples are represented by 
\[
\bz=\mathcal{S}\pi_{A,s}(\bw_{\text{aug}}+\widetilde{\bw_{\text{aug}}}).
\]
In this case, $\bw_{\text{aug}}^*=\bw_{\text{aug}}+\widetilde{\bw_{\text{aug}}}.$ Therefore
\begin{eqnarray*}
    \|\bw_{\text{aug}}^*-\bw_{\text{aug}}\|_2=\|\widetilde{\bw_{\text{aug}}}\|_2
    &\ge& 
    \frac{1}{\sqrt{(1+\delta)C_{A,k,s}}}\|\mathcal{W}\mathcal{P}\mathcal{S}\pi_{A,s}(\widetilde{\bw_{\text{aug}}})\|_2\\
    &=&
    \frac{1}{\sqrt{(1+\delta)C_{A,k,s}}}\|\mathcal{W}\mathcal{P}\widetilde{\be}\|_2.
\end{eqnarray*}
\end{enumerate}
\end{proof}

\begin{proof}[The proof of \Cref{thm:alg2}]
Let $\bw_{\text{aug}}^*$ be the solution of \eqref{recons2}, then for all $\overline{\bv}\in\mathbb{R}^{2n}$,  we have
\begin{eqnarray*}
&\quad& \|\mathcal{W}\mathcal{P}(\mathcal{S}\pi_{A,s}(\bw_{\text{aug}}^*)-\bz)\|^2_2+\gamma{\bw_{\text{aug}}^*}^{\top}g(\widetilde{L})\bw_{\text{aug}}^*\\
&\le& \|\mathcal{W}\mathcal{P}(\mathcal{S}\pi_{A,s}(\overline{\bv})-\bz)\|^2_2+\gamma{\overline{\bv}}^{\top}g(\widetilde{L})\overline{\bv}.  
\end{eqnarray*}

Notice that $\bw_{\text{aug}}^*=\alpha^*+\beta^*$ with $\alpha^*\in \text{span}({\widetilde{U}_k})$ and $\beta^*\in \text{span}({\widetilde{U}_k})^\perp$ where $\text{span}({\widetilde{U}_k})^\perp$ is the orthogonal complement of $\text{span}({\widetilde{U}_k}).$ Let $\widetilde{\overline{U}_k}=\begin{bmatrix}
    \overline{U}_k &\textbf{0}\\
    \textbf{0} & \overline{U}_k
\end{bmatrix}$ with $\overline{U}_k=(\textbf{u}_{k+1},\cdots,\textbf{u}_{n}),$ then $\text{span}({\widetilde{U}_k})^\perp=\text{span}(\widetilde{\overline{U}_k}).$ In addition, define $$g(\widetilde{\Theta_k})=\text{diag}(g(\theta_1),\cdots,g(\theta_k),g(\theta_1),\cdots,g(\theta_k))$$
and
$$g(\widetilde{\Theta_k}^\perp)=\text{diag}(g(\theta_{k+1}),\cdots,g(\theta_n),g(\theta_{k+1}),\cdots,g(\theta_n)).$$

In particular, if $\overline{\bv}= {\bw_{\text{aug}}}\in\text{span}(\widetilde{U}_k),$ we obtain
\begin{eqnarray*}
    &\quad& \|\mathcal{W}\mathcal{P}(\mathcal{S}\pi_{A,s}(\bw_{\text{aug}}^*)-\bz)\|^2_2+\gamma{\bw_{\text{aug}}^*}^{\top}g(\widetilde{L})\bw_{\text{aug}}^*\\
    &=& \|\mathcal{W}\mathcal{P}(\mathcal{S}\pi_{A,s}(\bw_{\text{aug}}^*)-\bz)\|^2_2+\gamma(\widetilde{U}_k^{\top}\alpha^*)^\top g(\widetilde{\Theta_k})(\widetilde{U}_k^{\top}\alpha^*)+\gamma(\widetilde{\overline{U}_k}^{\top}\beta^*)^\top g(\widetilde{\Theta_k}^\perp)(\widetilde{\overline{U}_k}^{\top}\beta^*)\\
&\le& \|\mathcal{W}\mathcal{P}(\mathcal{S}\pi_{A,s}(\bw_{\text{aug}})-\bz)\|^2_2+\gamma{\bw_{\text{aug}}}^{\top}g(\widetilde{L})\bw_{\text{aug}}\\
&=& \|\mathcal{W}\mathcal{P}\be\|^2_2+\gamma(\widetilde{U}_k^{\top}\bw_{\text{aug}})^\top g(\widetilde{\Theta_k})(\widetilde{U}_k^{\top}\bw_{\text{aug}})
\end{eqnarray*}
where we use the eigen-decompsition of $\widetilde{L}= \widetilde{U}_k \widetilde{\Theta_k} \widetilde{U}_k^\top + \widetilde{\overline{U}_k} \widetilde{\Theta_k}^\perp \widetilde{\overline{U}_k}^\top$ in the equalities.

Since $\theta_1\le\cdots\le\theta_n$,  the following inequalities hold: 
\[
(\widetilde{U}_k^{\top}\alpha^*)^\top g(\widetilde{\Theta_k})(\widetilde{U}_k^{\top}\alpha^*)\ge 0,
\]
\[
(\widetilde{\overline{U}_k}^{\top}\beta^*)^\top g(\widetilde{\Theta_k}^\perp)(\widetilde{\overline{U}_k}^{\top}\beta^*)\ge g(\theta_{k+1})\|\beta^*\|_2^2,
\]
and
\[
(\widetilde{U}_k^{\top}\bw_{\text{aug}})^\top g(\widetilde{\Theta_k})(\widetilde{U}_k^{\top}\bw_{\text{aug}})\le g(\theta_k)\|\bw_{\text{aug}}\|_2^2.
\]
Combining all the inequalities above, we have 
\begin{eqnarray}\label{a}
    \|\mathcal{W}\mathcal{P}(\mathcal{S}\pi_{A,s}(\bw_{\text{aug}}^*)-\bz)\|^2_2+\gamma g(\theta_{k+1})\|\beta^*\|_2^2\le \|\mathcal{W}\mathcal{P}\be\|^2_2+\gamma g(\theta_k)\|\bw_{\text{aug}}\|_2^2.
\end{eqnarray}
Since the left-hand side of $\eqref{a}$ is a sum of two positive quantities, it follows that
\begin{equation}
    \|\mathcal{W}\mathcal{P}(\mathcal{S}\pi_{A,s}(\bw_{\text{aug}}^*)-\bz)\|^2_2\le(\|\mathcal{W}\mathcal{P}\be\|_2+\sqrt{\gamma g(\theta_k)}\|\bw_{\text{aug}}\|_2)^2
\end{equation}
and
\begin{equation}
\gamma g(\theta_{k+1})\|\beta^*\|_2^2\le(\|\mathcal{W}\mathcal{P}\be\|_2+\sqrt{\gamma g(\theta_k)}\|\bw_{\text{aug}}\|_2)^2.
\end{equation}

Thus, we obtain
\begin{equation}\label{b}
    \|\mathcal{W}\mathcal{P}(\mathcal{S}\pi_{A,s}(\bw_{\text{aug}}^*)-\bz)\|_2\le\|\mathcal{W}\mathcal{P}\be\|_2+\sqrt{\gamma g(\theta_k)}\|\bw_{\text{aug}}\|_2
\end{equation}
and
\begin{equation}\label{c}
\sqrt{\gamma g(\theta_{k+1})}\|\beta^*\|_2\le\|\mathcal{W}\mathcal{P}\be\|_2+\sqrt{\gamma g(\theta_k)}\|\bw_{\text{aug}}\|_2.
\end{equation}

Furthermore, the following inequality holds
\begin{equation}\label{alpha}
\begin{aligned}
  \quad&\|\mathcal{W}\mathcal{P}(\mathcal{S}\pi_{A,s}(\bw_{\text{aug}}^*)-\bz)\|_2\\
  =&\|\mathcal{W}\mathcal{P}\mathcal{S}\pi_{A,s}(\bw_{\text{aug}}^*-\bw_{\text{aug}})-\mathcal{W}\mathcal{P}\be\|_2\\
=&\|\mathcal{W}\mathcal{P}\mathcal{S}\pi_{A,s}(\alpha^*+\beta^*-\bw_{\text{aug}})-\mathcal{W}\mathcal{P}\be\|_2\\
\ge&
\|\mathcal{W}\mathcal{P}\mathcal{S}\pi_{A,s}(\alpha^*-\bw_{\text{aug}})\|_2-\|\mathcal{W}\mathcal{P}\mathcal{S}\pi_{A,s}\beta^*\|_2-\|\mathcal{W}\mathcal{P}\be\|_2\\
\ge&\sqrt{(1-\delta)c_{A,k,s}}\|\alpha^*-\bw_{\text{aug}}\|_2-RC_{A,n,s}\|\beta^*\|_2-\|\mathcal{W}\mathcal{P}\be\|_2  
\end{aligned}
\end{equation}
where $\|\mathcal{W}\mathcal{P}\mathcal{S}\|_2\le R$ and $C_{A,n,s}$ is the upper bound of the embedding operator $\pi_{A,s}.$
Combining \eqref{b}, \eqref{c} and \eqref{alpha} gives
\begin{eqnarray*}  &\quad&\|\mathcal{W}\mathcal{P}\be\|_2+\sqrt{\gamma g(\theta_k)}\|\bw_{\text{aug}}\|_2\\
&\ge&
\sqrt{(1-\delta)c_{A,k,s}}\|\alpha^*-\bw_{\text{aug}}\|_2-RC_{A,n,s}\|\beta^*\|_2-\|\mathcal{W}\mathcal{P}\be\|_2\\
&\ge&\sqrt{(1-\delta)c_{A,k,s}}\|\alpha^*-\bw_{\text{aug}}\|_2-\|\mathcal{W}\mathcal{P}\be\|_2-RC_{A,n,s}\left(\frac{\|\mathcal{W}\mathcal{P}\be\|_2}{\sqrt{\gamma g(\theta_{k+1})}}+\sqrt{\frac{g(\theta_{k})}{g(\theta_{k+1})}}\|\bw_{\text{aug}}\|_2\right).
\end{eqnarray*}
Rearranging the above inequality, we obtain
\begin{eqnarray*}
    \|\alpha^*-\bw_{\text{aug}}\|_2\le\left(\frac{2+\frac{RC_{A,n,s}}{\sqrt{\gamma g(\theta_{k+1})}}}{\sqrt{(1-\delta)c_{A,k,s}}}\right)\|\mathcal{W}\mathcal{P}\be\|_2
    +\left(\frac{RC_{A,n,s}\sqrt{\frac{g(\theta_{k})}{g(\theta_{k+1})}}+\sqrt{\gamma g(\theta_k)}}{\sqrt{(1-\delta)c_{A,k,s}}}\right)\|\bw_{\text{aug}}\|_2.
\end{eqnarray*}
\end{proof}
\bibliographystyle{plain} 
\bibliography{ref}
\end{document}